\documentclass[12pt]{amsart}
\usepackage{amsfonts}

\usepackage{graphicx}
\usepackage{amsmath}
\usepackage{amsmath}
\usepackage{amsxtra}

\theoremstyle{plain}
\newtheorem{theorem}{Theorem}[section]
\newtheorem{lemma}[theorem]{Lemma}

\newtheorem{proposition}[theorem]{Proposition}
\newtheorem{corollary}[theorem]{Corollary}

\newtheorem{question}[theorem]{Question}
\theoremstyle{definition}

\newtheorem{example}[theorem]{Example}
\newtheorem{remark}[theorem]{Remark}

\numberwithin{equation}{section}
\begin{document}
\title{The extremal truncated moment problem}
\author{Ra\'{u}l E. Curto}
\address{Department of Mathematics\\
The University of Iowa\\
Iowa City, IA 52242-1419\\
USA}
\email{rcurto@math.uiowa.edu}
\author{Lawrence A. Fialkow}
\address{Department of Computer Science\\
State University of New York\\
New Paltz, NY 12561\\
USA}
\email{fialkowl@newpaltz.edu}
\author{H. Michael M\"{o}ller}
\address{FB Mathematik der Universit\"{a}t Dortmund\\
44221 Dortmund\\
Germany}
\email{michael.moeller@math.uni-dortmund.de}
\thanks{The first-named author's research was partially supported by NSF
Research Grants DMS-0099357 and DMS-0400741. \ The second-named author's
research was partially supported by NSF Research Grant DMS-0201430 and
DMS-0457138. }

\begin{abstract}
For a degree $2n$ real $d$-dimensional multisequence $\beta \equiv \beta
^{(2n)}$\newline
$=\{\beta _{i}\}_{i\in Z_{+}^{d},|i|\leq 2n}$ to have a \textit{representing
measure} $\mu $, it is necessary for the associated moment matrix ${\mathcal{%
M}}(n)(\beta )$ to be positive semidefinite and for the algebraic variety
associated to $\beta $, ${\mathcal{V}}\equiv {\mathcal{V}}_{\beta }$, to
satisfy $\operatorname{rank}\;{\mathcal{M}}(n)\leq \operatorname{card}\;{\mathcal{V}}$
as well as the following \textit{consistency} condition: if a polynomial $%
p(x)\equiv \sum_{|i|\leq 2n}a_{i}x^{i}$ vanishes on ${\mathcal{V}}$, then $%
\sum_{|i|\leq 2n}a_{i}\beta _{i}=0$. We prove that for the \textit{\ extremal%
} case ($\operatorname{rank}\;{\mathcal{M}}(n)=\operatorname{card}\;{\mathcal{V}}$),
positivity of ${\mathcal{M}}(n)$ and consistency are sufficient for the
existence of a (unique, $\operatorname{rank}\;{\mathcal{M}}(n)$-atomic)
representing measure. \ We also show that in the preceding result,
consistency cannot always be replaced by recursiveness of $\mathcal{M}(n)$.
\end{abstract}

\maketitle

\section{\label{sect1}Introduction}

Let $\beta \equiv \beta ^{(2n)}=\{\beta _{i}\}_{i\in \mathbb{Z}%
_{+}^{d},|i|\leq 2n}$ denote a real $d$-dimensional multisequence of degree $%
2n$. \ The \textit{truncated moment problem} for $\beta $ concerns the
existence of a positive Borel measure $\mu $, supported in $\mathbb{R}^{d}$,
such that 
\begin{equation}
\beta _{i}=\int_{\mathbb{R}^{d}}x^{i}~d\mu ,~~~|i|\leq 2n;  \label{eq1}
\end{equation}%
(here, for $x\equiv (x_{1},...,x_{d})\in \mathbb{R}^{d}$ and $i\equiv
(i_{1},...,i_{d})\in \mathbb{Z}_{+}^{d}$, we let $x^{i}:=x_{1}^{i_{1}}\cdots
x_{d}^{i_{d}})$. \ A measure $\mu $ as in (\ref{eq1}) is a \textit{%
representing measure} for $\beta $.

Let $\mathcal{P}\equiv \mathbb{R}^{d}[x]=\mathbb{R}[x_{1},...,x_{d}]$ denote
the space of real valued $d$-variable polynomials, and for $k\geq 1$, let $%
\mathcal{P}_{k}\equiv \mathbb{R}_{k}^{d}[x]$ denote the subspace of $%
\mathcal{P}$ consisting of polynomials $p$ with $\deg p\leq k$.
Corresponding to $\beta $ we have the \textit{Riesz functional} $\Lambda
\equiv \Lambda _{\beta }:\mathcal{P}_{2n}\rightarrow \mathbb{R}$, which
associates to an element $p$ of $\mathcal{P}_{2n}$, $p(x)\equiv
\sum_{|i|\leq 2n}a_{i}x^{i}$, the value $\Lambda (p):=\sum_{|i|\leq
2n}a_{i}\beta _{i}$; of course, in the presence of a representing measure $%
\mu $, we have $\Lambda (p)=\int p~d\mu $. \ In the sequel, $\hat{p}$
denotes the coefficient vector $(a_{i})$ of $p$.

Following \cite{tcmp1}, we associate to $\beta $ the \textit{moment matrix} $%
\mathcal{M}(n)\equiv \mathcal{M}(n)(\beta )$, with rows and columns $X^{i}$
indexed by the monomials of $\mathcal{P}_{n}$ in degree-lexicographic order;
for example, with $d=n=2$, the columns of $\mathcal{M}(2)$ are denoted as $%
1,X_{1},X_{2},X_{1}^{2},X_{2}X_{1},$ $X_{2}^{2}$. \ The entry in row $X^{i}$%
, column $X^{j}$ of $\mathcal{M}(n)$ is $\beta _{i+j}$, so $\mathcal{M}(n)$
is a real symmetric matrix characterized by 
\begin{equation}
\langle \mathcal{M}(n)\hat{p},\hat{q}\rangle =\Lambda (pq)~~~(p,q\in 
\mathcal{P}_{n}).  \label{eq2}
\end{equation}%
If $\mu $ is a representing measure for $\beta $, then $\langle \mathcal{M}%
(n)\hat{p},\hat{p}\rangle =\Lambda (p^{2})=\int p^{2}d\mu \geq 0$; since $%
\mathcal{M}(n)$ is real symmetric, it follows that $\mathcal{M}(n)$ is
positive semidefinite (in symbols, $\mathcal{M}(n)\geq 0)$.

The \textit{algebraic variety} of $\beta $ (or of $\mathcal{M}(n)(\beta )$)
is defined by 
\begin{equation*}
\mathcal{V}\equiv \mathcal{V}_{\beta }:=\bigcap_{p\in \mathcal{P}_{n},\hat{p}%
\in \ker \mathcal{M}(n)}~\mathcal{Z}(p),
\end{equation*}%
where $\mathcal{Z}(p):=\{x\in \mathbb{R}^{d}:p(x)=0\}$. \ (We sometimes
denote $\mathcal{V}_{\beta }$ as $\mathcal{V}(\mathcal{M}(n)(\beta )$.) \ If 
$\beta $ admits a representing measure $\mu $, then $p\in \mathcal{P}_{n}$
satisfies $\hat{p}\in \ker \mathcal{M}(n)$ if and only if $\operatorname{supp}%
\;\mu \subseteq \mathcal{Z}(p)$ \cite[Proposition 3.1]{tcmp1}. \ Thus $%
\operatorname{supp}\;\mu \subseteq \mathcal{V}$, and it follows from \cite[(1.7)]%
{tcmp3} that $r:=\operatorname{rank}\;\mathcal{M}(n)$ and $v:=\operatorname{card}\;%
\mathcal{V}$ satisfy $r\leq \operatorname{card}\;\operatorname{supp}\;\mu \leq v$.
Further, in this case, if $p\in \mathcal{P}_{2n}$ and $p|_{\mathcal{V}%
}\equiv 0$, then clearly $\Lambda (p)=\int p~d\mu =0$. \ To summarize the
preceding discussion, we have the following basic necessary conditions for
the existence of a representing measure for $\beta ^{(2n)}$: 
\begin{equation}
\text{(Positivity)}~~\mathcal{M}(n)\geq 0  \label{C1}
\end{equation}%
\begin{equation}
\text{(Consistency)}~~p\in \mathcal{P}_{2n},\;p|_{\mathcal{V}}\equiv
0\Longrightarrow \Lambda (p)=0  \label{C2}
\end{equation}%
\begin{equation}
\text{(Variety Condition)}~~r\leq v\text{, i.e., }\operatorname{rank}\;\mathcal{M}%
(n)\leq \operatorname{card}\;\mathcal{V}\text{.}  \label{C3}
\end{equation}%
As we show below (Section \ref{sect2}), consistency implies the following
condition: 
\begin{equation}
\text{(Recursiveness)}~~~p,q,pq\in \mathcal{P}_{n},\hat{p}\in \ker \mathcal{M%
}(n)\Longrightarrow \hat{pq}\in \ker \mathcal{M}(n).  \label{C4}
\end{equation}

Consistency is a new condition; previously, in \cite[p. 5]{tcmp1}, we
considered only recursiveness (when (\ref{C4}) holds, we say that $\beta $
(or $\mathcal{M}(n)(\beta )$) is \textit{recursively generated}). \ In %
\cite[Theorem 3.19]{tcmp1} we showed that for $d=1$ (the truncated \textit{%
Hamburger moment problem} for $\mathbb{R}$), positivity and recursiveness
are sufficient to imply the existence of representing measures. \ For $d=2$
(the plane), there exists $\mathcal{M}(3)>0$ (positive definite) for which $%
\beta $ has no representing measure \cite[Section 4]{tcmp2}. \ Since an
invertible moment matrix satisfies (\ref{C2}) and (\ref{C3}) vacuously, it
follows that in general (\ref{C1})-(\ref{C3}) are not sufficient conditions
for representing measures. \ By contrast, the results of \cite{tcmp5}, \cite%
{tcmp7}, and \cite{tcmp9} together show that when $d=2$ and $\ker \;\mathcal{%
M}(n)$ contains an element $\hat{p}$ with $\deg p\leq 2$, then $\beta $ has
a representing measure if and only if $\mathcal{M}(n)$ is positive,
recursively generated and satisfies the variety condition. \ This result
motivated the following question of \cite[Conjecture 1.2]{FiaIEOT}.

\begin{question}
\label{quest2}Suppose ${\mathcal{M}}(n)(\beta )$ is singular. If ${\mathcal{M%
}}(n)$ is positive, recursively generated, and $r\leq v$, does $\beta $
admit a representing measure?
\end{question}

In the present note we focus on the following refinement of Question \ref%
{quest2}.

\begin{question}
\label{quest1}Suppose ${\mathcal{M}}(n)(\beta )$ is singular. If ${\mathcal{M%
}}(n)$ is positive, $\beta $ is consistent, and $r\leq v$, does $\beta $
admit a representing measure?
\end{question}

Our main result provides an affirmative answer to Question \ref{quest1} in
the \textit{extremal} case, when $r=v$.

\begin{theorem}
\label{thm13}For $\beta \equiv \beta ^{(2n)}$ extremal, i.e., $r=v$, the
following are equivalent: \newline
(i)~ $\beta $ has a representing measure; \newline
(ii)~ $\beta $ has a unique representing measure, which is $\operatorname{rank}\;{%
\mathcal{M}}(n)$-atomic; \newline
(iii) ~ ${\mathcal{M}}(n)\geq 0$ and $\beta $ is consistent.
\end{theorem}

In many cases, the conditions of Theorem \ref{thm13} provide a concrete
solution to the extremal case of the truncated moment problem. \ Indeed,
only elementary linear algebra is required to verify that $\mathcal{M}(n)$
is positive semidefinite, to compute its rank, and to identify the
dependence relations which enter into the definition of the variety $%
\mathcal{V}$. Further, as we show in Section \ref{sect2}, if the points of
the variety can be computed exactly (which may be feasible in specific
examples by using computer algebra), then only elementary linear algebra is
required to verify that $\beta $ is consistent. \ The proof of Theorem \ref%
{thm13} is included in Theorem \ref{thm32} (Section \ref{sect3}), which also
provides a simple procedure for computing the unique representing measure
for $\beta $.

If the points of $\mathcal{V}(\mathcal{M}(n))$ are not known exactly, then
it may be difficult to verify consistency directly; for this reason, it is
of interest to identify cases in which recursiveness, which is easy to
check, actually implies consistency. \ In Sections \ref{sect2}, \ref{sect4}
and \ref{newsect6} we study the extent to which ``consistency'' in Theorem %
\ref{thm13} can be replaced by ``recursiveness,'' or by a simplified
consistency condition. \ Consider a planar moment matrix ${\mathcal{M}}%
(3)\geq 0$ with ${\mathcal{M}}(2)>0$ and a column dependence relation $%
Y=X^{3}$. \ In Section \ref{sect4} we show that if ${\mathcal{M}}(3)$ (as
above) is extremal with $r=v=7$, then recursiveness is indeed sufficient for
a representing measure. \ By contrast, in Section \ref{newsect6} we show
that for an extremal ${\mathcal{M}}(3)$ as above, but with $r=v=8$, it may
happen that there is no representing measure (Theorem \ref{thm62}). \ This
result provides a perhaps surprising negative answer to Question \ref{quest2}%
, and also shows that in general consistency is a strictly stronger property
than recursiveness. \ In Theorem \ref{thm63} we show that for the preceding $%
r=v=8$ extremal problem, consistency reduces to checking that $\Lambda (h)=0$
for a particular polynomial $h\in \mathbb{R}[x,y]$ of degree $4$.

We next observe that the extremal case is inherent in the truncated moment
problem. \ A recent result of C. Bayer and J. Teichmann \cite{BaTe}
(extending a classical theorem of V. Tchakaloff \cite{Tch} and its
successive generalizations in\ \cite{Mys}, \cite{tchaka} and \cite{tcmp8})
implies that if $\beta ^{(2n)}$ has a representing measure, then it has a
finitely atomic representing measure. \ In \cite{tcmp3} it was shown that $%
\beta ^{(2n)}$ has a finitely atomic representing measure if and only if $%
\mathcal{M}(n)$ admits an extension to a positive moment matrix $\mathcal{M}%
(n+k)$ (for some $k\geq 0$), which in turn admits a rank-preserving (i.e., 
\textit{flat}) moment matrix extension $\mathcal{M}(n+k+1)$. \ Further, %
\cite[Theorem 1.2]{tcmp10} shows that any flat extension $\mathcal{M}(n+k+1)$
is an extremal moment matrix for which there is a computable $\operatorname{rank}%
\;\mathcal{M}(n+k)$-atomic representing measure $\mu $. \ Clearly, $\mu $ is
also a finitely atomic representing measure for $\beta ^{(2n)}$, and every
finitely atomic representing measure for $\beta ^{(2n)}$ arises in this way.
\ In this sense, the existence \ of a representing measure for $\beta
^{(2n)} $ is intimately related to the solution of an extremal truncated
moment problem.

We conclude this section with two examples related to the extremal truncated
moment problem. \ In the first example we illustrate extremal truncated
moment problems of arbitrarily large degree. \ To ease the exposition of
this example, we will present it in terms of the \textit{truncated complex
moment problem}. \ Let $\gamma \equiv \gamma ^{(2n)}=\{\gamma
_{ij}\}_{i,j\in \mathbb{Z}_{+}^{d},|i|+|j|\leq 2n}$ denote a $d$-dimensional
complex multisequence of degree $2n$. \ The truncated complex moment problem
for $\gamma $ concerns the existence of a positive Borel measure $\nu $ on $%
\mathbb{C}^{d}$ such that 
\begin{equation}
\gamma _{ij}=\int_{\mathbf{C}^{d}}\bar{z}^{i}z^{j}d\nu ~~(i,j\in
Z_{+}^{d},|i|+|j|\leq 2n),  \label{eq17}
\end{equation}%
(where $z\equiv (z_{1},...,z_{d})$, $\bar{z}\equiv (\bar{z}_{1},...,\bar{z}%
_{d})\in \mathbb{C}^{d}$, $i\equiv (i_{1},...,i_{d})$, $j\equiv
(j_{1},...,j_{d})\in \mathbb{Z}_{+}^{d}$, and $\bar{z}^{i}z^{j}:=\bar{z}%
_{1}^{i_{1}}\cdots \bar{z}_{d}^{i_{d}}z_{1}^{j_{1}}\cdots z_{d}^{j_{d}}$). \
The Riesz functional for $\gamma $ is defined by $\Lambda _{\gamma }(\bar{z}%
^{i}z^{j}):=\gamma _{i,j}$. \ The mapping $\mathbb{C}^{d}\times \mathbb{C}%
^{d}\mapsto \mathbb{R}^{2d}\times \mathbb{R}^{2d}$ defined by $(z,\bar{z}%
)\mapsto (x,y)$ (where $x:=(z+\bar{z})/2$ and $y:=(z-\bar{z})/2i$) induces a
correspondence between truncated moment problems on $\mathbb{C}^{d}$ and
truncated moment problems on $\mathbb{R}^{2d}$. \ Under this correspondence, 
$\gamma $ is associated to a $2d$-dimensional real multisequence $\beta $
(also of degree $2n$) via the formula $\Lambda _{\beta
}((x,y)^{(k,j)}):=\Lambda _{\gamma }(((z+\bar{z})/2)^{k}((z-\bar{z}%
)/2i)^{j})\;\;(k,j\in \mathbb{Z}_{+}^{d},\left| k\right| +\left| j\right|
\leq 2n)$; we write $\beta \equiv \mathcal{S}(\gamma )$. \ Let $\mathbb{C}%
^{d}[z,\bar{z}]=\mathbb{C}[z_{1},...,z_{d},\bar{z}_{1},...,\bar{z}_{d}]$ and
let $\mathbb{C}_{k}^{d}[z,\bar{z}]$ denote the subspace of polynomials $p(z,%
\bar{z})$ with $deg~p\leq k$. \ The \textit{complex moment matrix} $%
M(n)\equiv M(n)(\gamma )$ has rows and columns indexed by monomials in $z$
and $\bar{z}$ up to degree $n$ in degree-lexicographic order, such that $%
\langle M(n)\hat{p},\hat{q}\rangle =\Lambda (p\bar{q})~~~(p,q\in \mathbb{C}%
_{n}^{d}[z,\bar{z}])$. The \textit{variety} of $\gamma $ is defined as $%
V(\gamma ):=\bigcap_{p\in \mathbb{C}_{n}^{d}[z,\bar{z}],\hat{p}\in
ker~M(n)}~Z(p)$, where $Z(p):=\{z\in \mathbb{C}^{d}:p(z,\bar{z})=0\}$. \ The
close connection between $M(n)(\gamma )$ and $\mathcal{M}(n)(\mathcal{S}%
(\gamma ))$ is described in detail in \cite[Section 2]{tcmp10}; in
particular, both moment matrices share the same positivity, rank,
recursiveness, and consistency, and, up to the identification $\mathbb{C}%
^{d}\approx \mathbb{R}^{2d}$, the same variety and representing measures. \
For this reason, results such as Theorem \ref{thm13} admit direct analogues
for the truncated complex moment problem. \ (For related instances of this,
the reader is referred to \cite[Theorems 2.19 and 2.21]{tcmp10}).

\begin{example}
\label{ex14}For $n>0$, we exhibit an extremal $\gamma \equiv \gamma ^{(2n)}$
in one complex variable with $\operatorname{rank}\;M(n)(\gamma )=\operatorname{card}%
\;V(\gamma )=2n$. \ The rows and columns of $M(n)$ are indexed by $1,Z,\bar{Z%
},...,Z^{n},$ $\bar{Z}Z^{n-1},...,\bar{Z}^{n-1}Z,\bar{Z}^{n}$. \ We set $%
\gamma _{ii}=1~(0\leq i\leq n)$, and for $0<a<1$, we set $\gamma
_{0,2n-1}=\gamma _{2n-1,0}:=a$ and $\gamma _{0,2n}=\gamma _{2n,0}:=1-a^{2}$;
the remaining $\gamma _{ij}$ equal $0$. \ For example, with $n=3$ we have 
\begin{equation*}
M(3)=\left( 
\begin{array}{cccccccccc}
1 & 0 & 0 & 0 & 1 & 0 & 0 & 0 & 0 & 0 \\ 
0 & 1 & 0 & 0 & 0 & 0 & 0 & 1 & 0 & 0 \\ 
0 & 0 & 1 & 0 & 0 & 0 & 0 & 0 & 1 & 0 \\ 
0 & 0 & 0 & 1 & 0 & 0 & 0 & 0 & 0 & a \\ 
1 & 0 & 0 & 0 & 1 & 0 & 0 & 0 & 0 & 0 \\ 
0 & 0 & 0 & 0 & 0 & 1 & a & 0 & 0 & 0 \\ 
0 & 0 & 0 & 0 & 0 & a & 1 & 0 & 0 & 1-a^{2} \\ 
0 & 1 & 0 & 0 & 0 & 0 & 0 & 1 & 0 & 0 \\ 
0 & 0 & 1 & 0 & 0 & 0 & 0 & 0 & 1 & 0 \\ 
0 & 0 & 0 & a & 0 & 0 & 1-a^{2} & 0 & 0 & 1%
\end{array}%
\right) .
\end{equation*}%
Observe that in the column space of $M(n)$ we have $\bar{Z}Z=1$, $\bar{Z}%
^{n}-Z^{n}=a(Z^{n-1}-\bar{Z}^{n-1})$, and a basis for the column space is
given by ${\mathcal{B}}\equiv \{1,Z,\bar{Z},Z^{2},\bar{Z}^{2},$ $...,$ $%
Z^{i},\bar{Z}^{i},...,Z^{n-1},$ $\bar{Z}^{n-1},Z^{n}\}$. It follows readily
that $M(n)$ is recursively generated. \ Note that $M_{{\mathcal{B}}}$, the
compression of $M(n)$ to the rows and columns indexed by ${\mathcal{B}}$, is
of the form $J\oplus \left( 
\begin{array}{cc}
1 & a \\ 
a & 1%
\end{array}%
\right) $, where $J$ is an identity matrix of size $2n-2$. \ Thus $M_{{%
\mathcal{B}}}$ is a positive definite matrix, with $\operatorname{rank}\;M_{{%
\mathcal{B}}}=2n$. \ Since $\operatorname{rank}\;M(n)=\operatorname{rank}\;M_{{%
\mathcal{B}}}$, it follows from \cite[Proposition 2.3]{RGWSI} that $M(n)$ is
positive semidefinite. \ (In the language of \cite{tcmp1}, $M(n)$ is a 
\textit{flat extension} of $M_{{\mathcal{B}}}$.) \ Now $\bar{Z}Z=1$, so we
may apply the analysis of the \textit{truncated trigonometric moment problem}
from \cite{tcmp5}. \ Since $M(n)$ is positive and recursively generated, $%
\bar{Z}Z=1$, and $\operatorname{rank}\;M(n)=2n$, \cite[Theorem 3.5]{tcmp5}
implies that $\gamma $ has a unique representing measure, which is $2n$%
-atomic; in particular, $\operatorname{card}\;V(\gamma )\geq \operatorname{rank}%
\;M(n)(\gamma )=2n$. \ Now $V(\gamma )$ consists of common solutions of the
equations $\bar{z}z=1$ and $\bar{z}^{n}-z^{n}=a(z^{n-1}-\bar{z}^{n-1})$, so $%
V(\gamma )\subseteq Z(p)$, where $p(z,\bar{z})=z^{2n}+az^{2n-1}-az-1$. \
Thus, $\operatorname{card}\;V(\gamma )\leq \operatorname{card}\;Z(p)\leq 2n$, and it
follows that $\operatorname{card}\;V(\gamma )=2n=\operatorname{rank}\;M(n)(\gamma )$,
whence $\gamma $ is extremal. \ \textbf{\qed}
\end{example}

The preceding example does not illustrate Theorem \ref{thm13}, because we
did not conclude that $\operatorname{card}\;V(\gamma )=\operatorname{rank}%
\;M(n)(\gamma )$ until after we had established the existence of a
representing measure using \cite{tcmp5}. \ Moment theory can sometimes be
used to estimate the number and location of the zeros of a prescribed
polynomial; indeed, as a by-product of Example \ref{ex14}, we see that the
polynomial $p(z)\equiv z^{2n}+az^{2n-1}-az-1\;\;(0<a<1)$ has $2n$ distinct
zeros, all in the unit circle. \ (In response to our question, Professor
Srdjan Petrovic has provided a direct proof of this fact.) \ 

The next example does illustrate how Theorem \ref{thm13} can be used to
solve an extremal problem; in particular, it shows how to verify consistency
and how to compute the unique representing measure.

\begin{example}
\label{ex15}Consider the $2$-dimensional real moment matrix 
\begin{equation*}
{\mathcal{M}}(2)=\left( 
\begin{array}{cccccc}
1 & 0 & 0 & 1/2 & 0 & 3/2 \\ 
0 & 1/2 & 0 & -5/4 & 0 & -3/4 \\ 
0 & 0 & 3/2 & 0 & -3/4 & 0 \\ 
1/2 & -5/4 & 0 & 45/8 & 0 & 3/8 \\ 
0 & 0 & -3/4 & 0 & 3/8 & 0 \\ 
3/2 & -3/4 & 0 & 3/8 & 0 & 45/8%
\end{array}%
\right) .
\end{equation*}%
We denote the rows and columns of ${\mathcal{M}}(2)$ as $%
1,~X,~Y,~X^{2},~YX,~Y^{2}$ and we denote the moment corresponding to $%
x^{i}y^{j}$ by $\beta _{ij}$. \ Since the upper left $4\times 4$ corner of ${%
\mathcal{M}}(2)$ is positive definite and we have column relations $%
YX=-(1/2)Y$ and $Y^{2}=2-4X-X^{2}$, it follows that ${\mathcal{M}}(2)$ is
positive semidefinite with $\operatorname{rank}\;{\mathcal{M}}(2)=4$. \ The
variety ${\mathcal{V}}\equiv {\mathcal{V}}_{\beta }$ consists of the common
zeros of $f(x):=yx+\frac{1}{2}y$ and $g(x):=y^{2}+x^{2}+4x-2$; these are the
points $w_{k}\equiv (x_{k},y_{k})$ $(1\leq k\leq 4)$, given by $x_{1}=x_{2}=-%
\frac{1}{2}$, $y_{1}=\frac{\sqrt{15}}{2}$, $y_{2}=-y_{1}$, $x_{3}=-2-\sqrt{6}
$, $x_{4}=-2+\sqrt{6}$, $y_{3}=y_{4}=0$, so $\beta ^{(4)}$ is extremal. \ We
next apply the method of Section \ref{sect2} to verify that $\beta $ is
consistent, and to this end we will compute a basis for ${\mathcal{I}}%
_{4}:=\{p\in {\mathcal{P}}_{4}:p|_{{\mathcal{V}}}\equiv 0\}$. \ Let $%
W_{4}\equiv W_{4}({\mathcal{V}})$ denote the matrix with 4 rows and 15
columns defined as follows. \ The columns are indexed by the monomials in ${%
\mathcal{P}}_{4}$ in degree-lexicographic order, and the entry in row $k$,
column $Y^{i}X^{j}$ is $y_{k}^{i}x_{k}^{j}$ $(1\leq k\leq 4,i,j\geq
0,i+j\leq 4)$. Clearly, a polynomial $p\equiv \sum_{0\leq i+j\leq
4}a_{ij}x^{i}y^{j}\in {\mathcal{P}}_{4}$ vanishes on ${\mathcal{V}}$ if and
only if $\hat{p}\equiv (a_{ij})\in ker~W_{4}$. \ Row-reducing $W_{4}$, we
obtain 
\begin{equation*}
W_{red}=\left( 
\begin{array}{ccccccccccccccc}
1 & 0 & 0 & 0 & 0 & 2 & 1 & 0 & -1 & 0 & -\frac{9}{2} & 0 & \frac{1}{2} & 0
& \frac{15}{2} \\ 
&  &  &  &  &  &  &  &  &  &  &  &  &  &  \\ 
0 & 1 & 0 & 0 & 0 & -4 & 0 & 0 & 2 & 0 & 1 & 0 & -1 & 0 & -15 \\ 
&  &  &  &  &  &  &  &  &  &  &  &  &  &  \\ 
0 & 0 & 1 & 0 & -\frac{1}{2} & 0 & 0 & \frac{1}{4} & 0 & \frac{15}{4} & 0 & -%
\frac{1}{8} & 0 & -\frac{15}{8} & 0 \\ 
&  &  &  &  &  &  &  &  &  &  &  &  &  &  \\ 
0 & 0 & 0 & 1 & 0 & -1 & -\frac{9}{2} & 0 & \frac{1}{2} & 0 & \frac{81}{4} & 
0 & -\frac{1}{4} & 0 & -\frac{15}{4} \\ 
&  &  &  &  &  &  &  &  &  &  &  &  &  & 
\end{array}%
\right) ,
\end{equation*}%
from which it follows that $\dim \ker W_{red}=11$. \ The form of $W_{red}$
implies that there is a basis for $\ker W_{red}$ $(=\ker W_{4})$ of the form 
$\{\hat{f}_{i}\}_{i=1}^{11}$, where $\hat{f}_{i}\equiv
(a_{i,1},...,a_{i,15}) $ satisfies $a_{i,4+j}=\delta _{ij}$ $(1\leq j\leq
11) $. \ By explicitly computing this basis, we derive the following basis
for ${\mathcal{I}}_{4}$: $f_{1}:=\frac{1}{2}y+yx$, $f_{2}:=-2+4x+x^{2}+y^{2}$%
, $f_{3}:=-1+\frac{9}{2}x^{2}+x^{3}$, $f_{4}:=-\frac{1}{4}y+yx^{2}$, $%
f_{5}:=1-2x-\frac{1}{2}x^{2}+y^{2}x$, $f_{6}:=-\frac{15}{4}y+y^{3}$, $f_{7}:=%
\frac{9}{2}-x-\frac{81}{4}x^{2}+x^{4}$, $f_{8}:=\frac{1}{8}y+yx^{3}$, $%
f_{9}:=-\frac{1}{2}+x+\frac{1}{4}x^{2}+y^{2}x^{2}$, $f_{10}:=\frac{15}{8}%
y+y^{3}x$, $f_{11}:=-\frac{15}{2}+15x+\frac{15}{4}x^{2}+y^{4}$. \ Using the
moment data, it is now straightforward to verify that $\Lambda _{\beta
}(f_{i})=0$ $(1\leq i\leq 11) $, so $\beta $ is consistent.\newline
\indent%
Theorem \ref{thm13} now implies that $\beta $ has a unique representing
measure. \ To compute this measure we follow the procedure described in the
proof of Theorem \ref{thm32}.\newline
Consider the following basis for the column space of $\mathcal{M}(2)$, ${%
\mathcal{B}}=\{1,~X,~Y,~X^{2}\}$. \ Let%
\begin{equation*}
V_{{\mathcal{B}}}\equiv V_{{\mathcal{B}}}[{\mathcal{V}}]:=\left( 
\begin{array}{cccc}
1 & 1 & 1 & 1 \\ 
x_{1} & x_{2} & x_{3} & x_{4} \\ 
y_{1} & y_{2} & y_{3} & y_{4} \\ 
x_{1}^{2} & x_{2}^{2} & x_{3}^{2} & x_{4}^{2}%
\end{array}%
\right)
\end{equation*}

We show in Lemma \ref{lem31} that $V_{{\mathcal{B}}}$ is necessarily
invertible, and in the proof of Theorem \ref{thm32} we show that the unique
representing measure for $\beta $ is of the form $\mu =\sum_{k=1}^{4}\rho
_{k}\delta _{w_{k}}$, whence the densities $\rho _{k}$ are determined by 
\begin{equation*}
(\rho _{1},\rho _{2},\rho _{3},\rho _{4})^{T}=V_{{\mathcal{B}}}^{-1}(\beta
_{00},\beta _{01},\beta _{10},\beta _{02})^{T}
\end{equation*}%
(where $(\cdot )^{T}$ denotes transpose). \ Using the given moment values,
we find $\rho _{1}=\rho _{2}=\frac{1}{5},~\rho _{3}=\frac{9\sqrt{2}-7\sqrt{3}%
}{30\sqrt{2}}\cong 0.0142262,~\rho _{4}=\frac{9\sqrt{2}+7\sqrt{3}}{30\sqrt{2}%
}\cong 0.585774$. \ \textbf{\qed}
\end{example}

\section{\label{sectideals}Real Ideals and Necessary Conditions}

If $\beta ^{(2n)}$ has a representing measure $\mu $, then the Riesz
functional 
\begin{equation*}
\Lambda \equiv \Lambda _{\beta }:\mathcal{P}_{2n}\rightarrow \mathbb{R}%
,\quad \Lambda (x^{i}):=\beta _{i}\;(=\int_{\mathbb{R}^{d}}x^{i}\;d\mu
\;\;(\left| i\right| \leq 2n)),
\end{equation*}%
is \textit{square positive}, that is, 
\begin{equation*}
p\in \mathcal{P}_{n}\Rightarrow \Lambda (p^{2})\geq 0
\end{equation*}%
(equivalently, $\mathcal{M}(n)(\beta )$ is positive semidefinite, cf. (\ref%
{eq2})). \ If we assume, in addition, that for a representing measure $\mu $
all moments 
\begin{equation*}
\int_{{\mathbb{R}}^{d}}x^{i}d\mu \;\;(i\in \mathbb{Z}_{+}^{d})
\end{equation*}%
are convergent, then we can extend $\Lambda $ to ${\mathcal{P}}$ by letting%
\begin{equation*}
\Lambda (x^{i}):=\int_{{\mathbb{R}}^{d}}x^{i}d\mu ,\quad i\in \mathbb{Z}%
_{+}^{d},
\end{equation*}%
thus obtaining a square positive functional over ${\mathcal{P}}$ (e.g., if $%
\mu $ is an $m$-atomic measure with support $\{w_{1},...,w_{m}\}\subseteq 
\mathbb{R}^{d}$, then $\Lambda (p)=\sum_{i=1}^{m}p(w_{i})\mu (\{w_{i}\})$
for all polynomials $p$). \ If $\Lambda _{\beta }$ does extend to a square
positive linear functional $\Lambda $ on $\mathcal{P}$, then, as shown in %
\cite{Moe}, the set 
\begin{equation*}
{\mathcal{I}}:=\{p\in \mathcal{P}:\Lambda (p^{2})=0\}
\end{equation*}%
is a \textit{real ideal}, i.e., it is an ideal ( $p_{1},p_{2}\in {\mathcal{I}%
}\Rightarrow p_{1}+p_{2}\in {\mathcal{I}}$ and $p\in {\mathcal{I}},\ q\in 
\mathcal{P}\Rightarrow pq\in {\mathcal{I}}$) and satisfies one of the
following two equivalent conditions:%
\begin{equation*}
\begin{array}{l}
\text{(i) \ For \ }s\in \mathbb{Z}_{+},p_{1},\ldots ,p_{s}\in \mathcal{P}%
:\sum_{i=1}^{s}p_{i}^{2}\in {\mathcal{I}}\Rightarrow \{p_{1},\ldots
,p_{s}\}\subseteq {\mathcal{I}}\text{;} \\ 
\text{(ii) \ There exists}\ G\subseteq \mathbb{R}^{d}\text{ such that for
all }p\in \mathcal{P}:p|_{G}\equiv 0\Rightarrow p\in {\mathcal{I}}.%
\end{array}%
\end{equation*}%
If ${\mathcal{I}}$ is a real ideal, then one may take for $G$ the \textit{%
real variety} 
\begin{equation*}
V_{\mathbb{R}}({\mathcal{I}}):=\{w\in \mathbb{R}^{d}:f(w)=0\quad (\text{all }%
f\in {\mathcal{I}}\mathcal{)}\}.
\end{equation*}%
But one may also take for $G$ any subset of $V_{\mathbb{R}}({\mathcal{I}})$
containing sufficiently many points, such that 
\begin{equation*}
p\in {\mathcal{P}},\ p|_{G}\equiv 0\Rightarrow p|_{V_{\mathbb{\mathbb{%
\mathbb{\mathbb{R}}}}}({\mathcal{I}})}\equiv 0.
\end{equation*}%
For instance, if the real variety is a (real) line, one may take for $G$ a
subset of infinitely many points on that line. \ On the other hand, if $V_{%
\mathbb{R}}({\mathcal{I}})$ is a finite set of points, then necessarily $%
G=V_{\mathbb{R}}({\mathcal{I}})$. \ (We note that in the full moment problem
for $\beta \equiv \beta ^{(\infty )}$, M. Laurent \cite{Lau2} independently
showed that $\mathcal{J}:=\{p\in \mathcal{P}:M(\infty )\hat{p}=0\}$ is a 
\textit{radical} ideal; equivalently, $p\in \mathcal{J}\Leftrightarrow
p^{2}\in \mathcal{J}$.)

If ${\mathcal{I}}$ is an ideal, its subset ${\mathcal{I}}_{k}:={\mathcal{I}}%
\cap \mathcal{P}_{k}$ is an $\mathbb{R}$-vector subspace of $\mathcal{P}_{k}$%
. \ One can then introduce the \textit{Hilbert function} of ${\mathcal{I}}$
by 
\begin{equation*}
H_{{\mathcal{I}}}(k):=\dim \mathcal{P}_{k}-\dim {\mathcal{I}}_{k},\quad k\in 
\mathbb{Z}_{+};
\end{equation*}%
in \cite{CLO'S} this is called the \textit{affine Hilbert function}. \ As
shown for instance in \cite{CLO'S}, both $k\mapsto \dim {\mathcal{I}}_{k}$
and $k\mapsto H_{{\mathcal{I}}}(k)$ are nondecreasing functions, and for
sufficiently large $k$, say $k\geq k_{0},$ $H_{{\mathcal{I}}}(k)$ becomes a
polynomial in $k$, the so-called \textit{Hilbert polynomial of} ${\mathcal{I}%
}$, whose degree equals the \textit{dimension} of ${\mathcal{I}}$.

\begin{example}
\label{ex21}Let $G\equiv \{w_{1},\ldots ,w_{m}\}\subseteq \mathbb{R}^{d}$. \
Then ${\mathcal{I}}:=\{f\in \mathcal{P}:f|_{G}\equiv 0\}$ is a real ideal
with $V_{\mathbb{R}}({\mathcal{I}})=G$. \ Let $t_{1},t_{2},t_{3},\ldots $
denote the monomials $x^{i}$ in degree-lexicographic order, so that for each 
$k\in \mathbb{Z}_{+}$ $t_{1},\ldots ,t_{K}$ (with $K:=\dim \mathcal{P}_{k}$)
form a basis of the $\mathbb{R}$-vector space $\mathcal{P}_{k}$. \ For $p\in 
\mathcal{P}_{k}$, $p\equiv \sum_{i=1}^{K}a_{i}t_{i}$, let $\hat{p}%
:=(a_{1},\ldots ,a_{K})$ (the coefficient vector of $p$). \ Then $p(x)$ can
be written as 
\begin{equation*}
p(x)=\left\langle \hat{p},t(x)\right\rangle ,
\end{equation*}%
where $t(x):=(t_{1}(x),...,t_{K}(x))$, so 
\begin{equation*}
p\in {\mathcal{I}}\cap \mathcal{P}_{k}\Leftrightarrow \hat{p}\perp
t(w_{i}),\ i=1,\ldots ,m.
\end{equation*}%
Arranging the rows $t(w_{i})\;(=(t_{1}(w_{i}),\ldots ,t_{K}(w_{i})))$ in a
matrix 
\begin{equation*}
W_{k}\equiv W_{k}[G]:=(t_{j}(w_{i}))_{i=1,\ldots ,m,\ j=1,\ldots ,K},
\end{equation*}%
one gets $p\in {\mathcal{I}}\cap \mathcal{P}_{k}\Leftrightarrow \hat{p}\in
\ker \;W_{k}$, whence $\dim {\mathcal{I}}_{k}+\operatorname{rank}\;W_{k}=\dim 
\mathcal{P}_{k}$, or using the Hilbert function, 
\begin{equation*}
H_{{\mathcal{I}}}(k)=\operatorname{rank}\;W_{k},\ \;k\in \mathbb{Z}_{+}\;.
\end{equation*}%
By construction, $W_{k}$ is a submatrix of $W_{k+1}$. \ Hence $\operatorname{rank}%
\;W_{k}\leq \operatorname{rank}\;W_{k+1}$, reflecting the fact that the Hilbert
function increases. \ If, for a given $k$, the rank of $W_{k}$ is less than $%
m$, then one row of $W_{k}$, say the last one, depends on the others. \ This
means that every polynomial which vanishes in $w_{1},\ldots ,w_{m-1}$ also
vanishes in $w_{m}$. \ Using Lagrange interpolation polynomials, we see that
for all sufficiently large $k$ this cannot happen. \ Hence $\operatorname{rank}%
\;W_{k}=m$ for all sufficiently large $k$. \ This $m$ is the constant
(degree-$0$) polynomial in $k$ which coincides with $H_{{\mathcal{I}}}(k)$
for all $k\geq k_{0}$; hence, ${\mathcal{I}}$ is a zero dimensional ideal. \ 
\textbf{\qed}
\end{example}

Now we will study the consistency condition (\ref{C2}). \ We consider an
arbitrary real $d$-dimensional multisequence $\beta \equiv \beta ^{(2n)}$ of
degree $2n$. \ Associated with $\beta $ one has the Riesz functional $%
\Lambda $, the moment matrix $\mathcal{M}(n)$, and the algebraic variety $%
\mathcal{V}\equiv \mathcal{V}_{\beta }$ (or $\mathcal{V}(\mathcal{M}(n))$).
\ One can then define the ideal 
\begin{equation}
\mathcal{I}(\mathcal{V}):=\{p\in \mathcal{P}:p|_{\mathcal{V}}\equiv 0\}.
\label{iv}
\end{equation}%
Since $\mathcal{V}$ is a set of real points, $\mathcal{I}(\mathcal{V})$ is a
real ideal, which we will call the \textit{real ideal} of $\beta $. \ 

\begin{lemma}
\label{lem22}Assume that $\beta \equiv \beta ^{(2n)}$ satisfies (\ref{C2}).
\ Then 
\begin{equation}
\mathcal{N}_{n}:=\{p\in \mathcal{P}_{n}:\mathcal{M}(n)\hat{p}=0\}=\mathcal{I}%
(\mathcal{V})\bigcap \mathcal{P}_{n}.  \label{eq21}
\end{equation}%
If $t_{1},...,t_{N}$ denote the monomials $x^{i}\in \mathcal{P}_{n}$ in
degree-lexicographic order, then the row vectors of $\mathcal{M}(n)$ and the
row vectors $\{t(w):=(t_{1}(w),...,t_{N}(w)):w\in \mathcal{V}\}$, span the
same subspace of $\mathbb{R}^{N}$; in particular, $\operatorname{rank}\;\mathcal{M%
}(n)=H_{\mathcal{I}(\mathcal{V})}(n)$.
\end{lemma}

\begin{proof}
If $p\in \mathcal{I}(\mathcal{V})\bigcap \mathcal{P}_{n}$ and $q\in \mathcal{%
P}_{n}$, then $pq\in \mathcal{P}_{2n}$ and $(pq)|_{\mathcal{V}}\equiv 0$,
whence by the consistency property (\ref{C2}) we must have $\left\langle 
\mathcal{M}(n)\hat{p},\hat{q}\right\rangle =\Lambda (pq)=0$; thus, $\mathcal{%
M}(n)\hat{p}=0$. \ Conversely, if $p\in \mathcal{P}_{n}$ and $\mathcal{M}(n)%
\hat{p}=0$, then $p|_{\mathcal{V}}\equiv 0$ by the definition of $\mathcal{V}
$. \ Hence $p\in \mathcal{I}(\mathcal{V})$. \ Now, using (\ref{eq21}) and
proceeding as in Example \ref{ex21}, we see that 
\begin{eqnarray*}
\hat{p} &\in &\ker \;\mathcal{M}(n)\Leftrightarrow p\in \mathcal{I}(\mathcal{%
V})\bigcap \mathcal{P}_{n} \\
&\Leftrightarrow &\hat{p}\perp t(w)\;\;(\text{all }w\in \mathcal{V}).
\end{eqnarray*}%
This means that the rows of $\mathcal{M}(n)$ span the same space (namely, $%
\mathbb{R}^{N}\ominus \ker \;\mathcal{M}(n)$) as the rows $(t_{1}(w),\ldots
,t_{N}(w))$, $w\in \mathcal{V}$. $\ $It also follows that 
\begin{equation*}
\operatorname{rank}\;{\mathcal{M}}(n)=\dim \;{\mathcal{P}}_{n}-\dim \;\ker \;{%
\mathcal{M}}(n)=\dim \;{\mathcal{P}}_{n}-\dim \;\mathcal{I}(\mathcal{V})\cap 
{\mathcal{P}}_{n}=H_{\mathcal{I}(\mathcal{V})}(n).
\end{equation*}
\end{proof}

As the following lemma will show, consistency is a very strong condition,
already yielding an atomic measure (though one which may have some negative
densities).

\begin{lemma}
\label{lem23}Let $\Lambda :\mathcal{P}_{2n}\rightarrow \mathbb{R}$ be a
linear functional and let $\mathcal{V}\subseteq \mathbb{R}^{d}$. \ The
following statements are equivalent.\newline
(a) There exist $\alpha _{1},...,\alpha _{m}\in \mathbb{R}$ and there exist $%
w_{1},...,w_{m}\in \mathcal{V}$ such that $\Lambda (p)=\sum_{i=1}^{m}\alpha
_{i}p(w_{i})\;($all $p\in \mathcal{P}_{2n}).$\newline
(b) If $p\in \mathcal{P}_{2n}$ and $p|_{\mathcal{V}}\equiv 0$, then $\Lambda
(p)=0$.
\end{lemma}

\begin{proof}
The implication (a) ${\Rightarrow }$ (b) is obvious. \ Therefore assume that
(b) holds, and fix the basis of monomials $x^{i}$ of $\mathcal{P}_{2n}$. \
For notational convenience, denote this basis by $t_{1},\ldots ,t_{K}$. \
Then b) is equivalent to 
\begin{equation*}
\text{(c) \ For all }c_{1},...,c_{K}\in \mathbb{R}^{K}:%
\sum_{j=1}^{K}c_{j}t_{j}(w)=0\;\;(\text{all }w\in \mathcal{V})\Rightarrow
\sum_{j=1}^{K}c_{j}\Lambda (t_{j})=0.
\end{equation*}%
Using $\hat{c}:=(c_{1},\ldots ,c_{K}),\ t(w):=(t_{1}(w),\ldots ,t_{K}(w)),$
and $\hat{\Lambda}:=(\Lambda (t_{1}),\ldots ,\Lambda (t_{K}))$, (b) is thus
equivalent to 
\begin{equation*}
\text{(d) \ For all }\hat{c}\in \mathbb{R}^{K}:\hat{c}\perp t(w)\;\;(\text{%
all }w\in \mathcal{V})\Rightarrow \hat{c}\perp \hat{\Lambda}.
\end{equation*}%
Recall that for subspaces $\mathcal{R}$ and $\mathcal{S}$ of $\mathbb{R}^{K}$%
, $\mathcal{R}^{\perp }\subseteq \mathcal{S}^{\perp }\Leftrightarrow 
\mathcal{S}\subseteq \mathcal{R}$. \ Hence $\hat{\Lambda}$ is in the $%
\mathbb{R}$-linear subspace of $\mathbb{R}^{K}$ spanned by $\{t(w):w\in 
\mathcal{V}\}$. \ As such, this subspace has a basis of $m\;(\leq K)$
vectors $t(w_{1}),...,t(w_{m})$, where $w_{1},...,w_{m}\in \mathcal{V}$. \
Hence there exist $\alpha _{1},...,\alpha _{m}\in \mathbb{R}$ such that $%
\hat{\Lambda}=\sum_{i=1}^{m}\alpha _{i}t(w_{i})$, or equivalently, 
\begin{equation*}
\Lambda (t_{j})=\sum_{i=1}^{m}\alpha _{i}t_{j}(w_{i})\;\;(1\leq j\leq K).
\end{equation*}%
This is a linear relation holding for a basis of $\mathcal{P}_{2n}$, hence
it holds true for all $p\in \mathcal{P}_{2n}$, that is,%
\begin{equation*}
p\in \mathcal{P}_{2n}\Rightarrow \Lambda (p)=\sum_{i=1}^{m}\alpha
_{i}p(w_{i}).
\end{equation*}
\end{proof}

\begin{remark}
\label{rem1}If $\Lambda $ is the Riesz functional $\Lambda _{\beta }$
corresponding to $\beta \equiv \beta ^{(2n)}$, then Lemma \ref{lem23}(b) is
the consistency condition (\ref{C2}). \ We remark that in the proof of Lemma %
\ref{lem23}(b) we did not assume the square positivity of $\Lambda $ (which
corresponds to the positivity condition (\ref{C1}) when $\Lambda =\Lambda
_{\beta }$).\newline
\indent%
When $\Lambda =\Lambda _{\beta }$, $\mathcal{V}\equiv \mathcal{V}(\mathcal{M}%
(n))=\{w_{1},...,w_{m}\}$, and $\operatorname{rank}\;\mathcal{M}(n)=m$ (the
extremal case), we next show that in the representation of Lemma \ref{lem23}%
(a), the square positivity of $\Lambda $ is equivalent to the positivity of
the $\alpha _{i}$'s. \ We have noted above that $\Lambda $ is square
positive if and only if $\mathcal{M}(n)$ is positive semidefinite; in this
case, we also have $\{p\in \mathcal{P}_{n}:\mathcal{M}(n)\hat{p}=0\}=\{p\in 
\mathcal{P}_{n}:\Lambda (p^{2})=0\}$.
\end{remark}

\begin{lemma}
\label{lem24}Let $\Lambda \equiv \Lambda _{\beta }:\mathcal{P}%
_{2n}\rightarrow \mathbb{R}$ be given by 
\begin{equation*}
\Lambda (p):=\sum_{i=1}^{m}\alpha _{i}p(w_{i})\;\;(p\in \mathcal{P}_{2n}),
\end{equation*}%
with $\mathcal{V}_{\beta }\equiv \{w_{1},...,w_{m}\}\subseteq \mathbb{R}^{d}$
and $\alpha _{1},...,\alpha _{m}\in \mathbb{R}$. \ If $\operatorname{rank}\;%
\mathcal{M}(n)=m$, the following statements are equivalent:\newline
(i) \ $\alpha _{i}>0\;\;($all $i=1,...,m)$;\newline
(ii) \ $\Lambda $ is square positive.
\end{lemma}

\begin{proof}
The implication $(i)\Rightarrow (ii)$ is obvious. \ Conversely, assume that $%
\Lambda $ is square positive, i.e., $\mathcal{M}(n)$ is positive
semi-definite. \ Let $t_{1},...,t_{N}$ be the basis of monomials in $%
\mathcal{P}_{n}$ in degree-lexicographic order, so that the $(j,k)$-entry of 
$\mathcal{M}(n)$ is $\Lambda (t_{j}t_{k})$. \ It follows that $\mathcal{M}%
(n) $ can be decomposed as 
\begin{equation}
\mathcal{M}(n)=W_{n}^{T}\left( 
\begin{array}{ccc}
\alpha _{1} &  & 0 \\ 
& \ddots &  \\ 
0 &  & \alpha _{m}%
\end{array}%
\right) W_{n},  \label{new23}
\end{equation}%
where $W_{n}$ is the $m\times N$ matrix with rows $t(w_{i})\equiv
(t_{1}(w_{i}),...,t_{N}(w_{i}))\;\;(1\leq i\leq m)$. \ Since $\operatorname{rank}%
\;\mathcal{M}(n)$ $=m$, (\ref{new23}) implies that $\operatorname{rank}\;W_{n}=m$%
. \ Hence the columns of $W_{n}$ span $\mathbb{R}^{m}$; in particular, every
unit vector in $\mathbb{R}^{m}$ is a linear combination of columns of $W_{n}$%
. \ This means that there exist polynomials $\ell _{i}\in \mathcal{P}_{n}$
satisfying $\ell _{i}(w_{j})=\delta _{ij}\;\;(1\leq i,j\leq m)$, where $%
\delta _{ij}$ denotes the Kronecker symbol. \ Now, $\alpha _{i}=\Lambda
(\ell _{i}^{2})=\left\langle \mathcal{M}(n)\hat{\ell}_{i},\hat{\ell}%
_{i}\right\rangle \geq 0$ (since $\mathcal{M}(n)\geq 0$). \ Finally, no $%
\alpha _{i}$ can be zero, because otherwise $\operatorname{rank}\;\mathcal{M}%
(n)<m $, a contradiction.
\end{proof}

\begin{remark}
\label{rem2}(i) The preceding results yield a first proof of Theorem \ref%
{thm13}(iii) $\Rightarrow $ (i). \ Indeed, Lemma \ref{lem23} shows that if $%
\beta $ is consistent, then $\beta $ admits an atomic representing measure $%
\mu $, while Lemma \ref{lem24} shows that if $\mathcal{M}(n)$ is also
positive semi-definite and extremal, then $\mu $ is $\operatorname{rank}\;%
\mathcal{M}(n)$-atomic and $\mu \geq 0$.\newline
(ii) A decomposition similar to (\ref{new23}) was used by Laurent \cite{Lau2}
in her study of the full moment problem for $\beta ^{(\infty )}$ in the case
when $\operatorname{card}\;\mathcal{V}(\mathcal{M}(\infty ))<+\infty $.
\end{remark}

We conclude this section with some additional observations about ideals and
consistency. \ Given a real $d$-dimensional multisequence $\beta $ of degree 
$2n$, let $\{p_{1},...,p_{s}\}$ denote a basis for $\mathcal{N}_{n}:=\{p\in 
\mathcal{P}_{n}:\mathcal{M}(n)\hat{p}=0\}$. \ Denote by $\mathcal{J}\equiv 
\mathcal{J}_{\beta }$ the smallest ideal containing the polynomials $%
p_{1},...,p_{s}$. \ Since $\mathcal{V}\equiv \mathcal{V}_{\beta }$ is the
set of all real common zeros of $p_{1},...,p_{s}$, we have $\mathcal{J}%
\subseteq \mathcal{I}(\mathcal{V})$. \ If $\beta $ is consistent, then Lemma %
\ref{lem22} gives $\mathcal{J}\bigcap \mathcal{P}_{n}=\mathcal{I}(\mathcal{V}%
)\bigcap \mathcal{P}_{n}$, whence 
\begin{equation*}
\dim (\mathcal{J}\bigcap \mathcal{P}_{k})\leq \dim (\mathcal{I}(\mathcal{V}%
)\bigcap \mathcal{P}_{k})\;\;(k\geq 0),
\end{equation*}%
with equality when $k=0,...,n$. \ For general $\beta $, the consistency
condition (\ref{C2}) can be rephrased in terms of $\mathcal{I}(\mathcal{V})$
as 
\begin{equation*}
p\in \mathcal{I}(\mathcal{V})\bigcap \mathcal{P}_{2n}\Rightarrow \Lambda
_{\beta }(p)=0.
\end{equation*}%
Now, since $\mathcal{J}\bigcap \mathcal{P}_{2n}$ is a subset of $\mathcal{I}(%
\mathcal{V})\bigcap \mathcal{P}_{2n}$, we can find 
\begin{equation*}
M:=\dim (\mathcal{I}(\mathcal{V})\bigcap \mathcal{P}_{2n})-\dim (\mathcal{J}%
\bigcap \mathcal{P}_{2n})\;(=H_{\mathcal{I}(\mathcal{V})}(2n)-H_{\mathcal{J}%
}(2n))
\end{equation*}%
polynomials $h_{1},...,h_{M}\in \mathcal{I}(\mathcal{V})\bigcap \mathcal{P}%
_{2n}$ enlarging a basis for $\mathcal{J}\bigcap \mathcal{P}_{2n}$ to a
basis for $\mathcal{I}(\mathcal{V})\bigcap \mathcal{P}_{2n}$. \ Then (\ref%
{C2}) can be rephrased again as 
\begin{equation}
p\in \mathcal{J}\bigcap \mathcal{P}_{2n}\Rightarrow \Lambda (p)=0\text{, and 
}\Lambda (h_{i})=0\;\;(1\leq i\leq M).  \label{new24}
\end{equation}%
Note that if $f\in \mathcal{N}_{n}$ and $g\in \mathcal{P}_{n}$, then $%
p:=fg\in \mathcal{J}\bigcap \mathcal{P}_{2n}$ and $\Lambda (p)=\left\langle 
\mathcal{M}(n)\hat{f},\hat{g}\right\rangle =0$. \ In Sections \ref{sect2}
and \ref{newsect6} we will identify situations in which $p\in \mathcal{J}%
\bigcap \mathcal{P}_{2n}$ always implies $\Lambda (p)=0$, so that
consistency reduces to the test $\Lambda (h_{i})=0\;\;(1\leq i\leq M)$.

\section{\label{sect2}Moment Matrices and Consistency}

A basic result of \cite{tcmp1} shows that $\beta \equiv \beta ^{(2n)}$ has a 
\textit{minimal representing measure}, i.e., a representing measure whose
support consists of exactly $\operatorname{rank}\;\mathcal{M}(n)$ atoms, if and
only if $\mathcal{M}(n)\geq 0$ and $\mathcal{M}(n)$ admits an extension to a
moment matrix $\mathcal{M}(n+1)$ with $\operatorname{rank}\;\mathcal{M}(n+1)=%
\operatorname{rank}\;\mathcal{M}(n)$. Following \cite{tcmp1}, we refer to such an
extension as a \textit{flat extension}. \ There is at present no concrete
set of necessary and sufficient conditions for the existence of flat
extensions $\mathcal{M}(n+1)$; one useful sufficient condition is that $%
\mathcal{M}(n)\geq 0$ satisfy $\operatorname{rank}\;\mathcal{M}(n)=\operatorname{rank}%
\;\mathcal{M}(n-1)$ \cite[Theorem 5.4]{tcmp1}. \ More generally, $\beta $
has a \textit{finitely atomic} representing measure (a representing measure
with finite support) if and only if $\mathcal{M}(n)$ admits a positive
extension $\mathcal{M}(n+k)$ (for some $k\geq 0$), which in turn admits a
flat extension $\mathcal{M}(n+k+1)$ (cf. \cite[Theorem 1.5]{tcmp3}). \ Since 
$\mathcal{M}(n+k+1)$ then admits unique successive flat extensions $\mathcal{%
M}(n+k+2)$, $\mathcal{M}(n+k+3)$, ... \cite{tcmp1}, this condition is
equivalent to the existence of a finite rank positive extension $\mathcal{M}%
(\infty )$. \ Further, a recent result of C. Bayer and J. Teichmann \cite%
{BaTe} (cf. Section \ref{sect1}) implies that if $\beta $ has a representing
measure, then $\beta $ admits a finitely atomic representing measure as just
described.

Recall that the columns of $\mathcal{M}(n)$ are denoted as $X^{i}$, $|i|\leq
n$, following the degree-lexicographic ordering of the monomials $x^{i}$ in $%
\mathcal{P}_{n}$. \ Let $p\in \mathcal{P}_{n}$, $p(x)\equiv \sum_{|i|\leq
n}a_{i}x^{i}$; the general element of $\mathcal{C}_{\mathcal{M}(n)}$, the
column space of $\mathcal{M}(n)$, may thus be denoted as $%
p(X):=\sum_{|i|\leq n}a_{i}X^{i}$. \ Let $\hat{p}\equiv (a_{i})$ denote the
coefficient vector of $p$ relative to the basis of monomials of $\mathcal{P}%
_{n}$ in degree-lexicographic order, and note that $p(X)=\mathcal{M}(n)\hat{p%
}$. \ Now recall the variety of $\beta $, 
\begin{equation*}
\mathcal{V}\equiv \mathcal{V}_{\beta }:=\bigcap_{p\in \mathcal{P}_{n},p(X)=0}%
\mathcal{Z}(p),
\end{equation*}%
where $\mathcal{Z}(p):=\{x\in \mathbb{R}^{d}:p(x)=0\}$. \ Let $\mathcal{P}%
_{n}|_{\mathcal{V}}$ denote the restriction to $\mathcal{V}$ of the
polynomials in $\mathcal{P}_{n}$, and consider the mapping $\phi _{\beta }:%
\mathcal{C}_{\mathcal{M}(n)}\rightarrow \mathcal{P}_{n}|_{\mathcal{V}}$
given by $p(X)\mapsto p|_{\mathcal{V}}$. \ The map $\phi _{\beta }$ is
well-defined, for if $p,q\in \mathcal{P}_{n}$ with $p(X)=q(X)$, then $%
\mathcal{V}\subseteq \mathcal{Z}(p-q)$, whence $p|_{\mathcal{V}}=q|_{%
\mathcal{V}}$. \ Note that if $\beta $ has a representing measure $\mu $,
then $\phi _{\beta }$ is 1-1; for, if $p\in \mathcal{P}_{n}$ and $p|_{%
\mathcal{V}}\equiv 0$, then since $\operatorname{supp}\;\mu \subseteq \mathcal{V}$
(cf. Section \ref{sect1}), we have $p|_{\operatorname{supp}\;\mu }\equiv 0$,
whence \cite[Proposition 3.1]{tcmp1} implies $p(X)=0$. \ Consider also the
following property of $\beta $:%
\begin{equation}
p\in \mathcal{P}_{n},q\in \mathcal{P},pq\in \mathcal{P}_{2n},p(X)=0%
\Rightarrow \Lambda _{\beta }(pq)=0  \label{neweq31}
\end{equation}%
(where $\Lambda _{\beta }$ is the Riesz functional associated to $\beta $;
cf. Section \ref{sect1}).

The following result will be used in the proof of Theorem \ref{thm13}.

\begin{proposition}
\label{consistency}Let $\beta ,\phi _{\beta }$ and $\mathcal{M}(n)(\beta )$
be as above. \ Then\newline
(i) \ $\beta $ consistent $\Longrightarrow \phi _{\beta }$ 1-1 $%
\Longrightarrow \mathcal{M}(n)(\beta )$ recursively generated.\newline
(ii) \ $\beta $ consistent $\Longrightarrow \beta $ satisfies (\ref{neweq31}%
) $\Longrightarrow \mathcal{M}(n)(\beta )$ recursively generated.
\end{proposition}

\begin{proof}
(i) \ Suppose $\beta $ is consistent. \ Formula (\ref{eq21}) in Lemma \ref%
{lem22} implies that $\phi _{\beta }$ is 1-1. \ We next assume that $\phi
_{\beta }$ is 1-1 and we show that $\mathcal{M}(n)$ is recursively
generated. \ Let $p,~q,~pq\in \mathcal{P}_{n}$ and suppose $p(X)=0$. \ Since 
$\mathcal{V}\subseteq \mathcal{Z}(p)$, then $p|_{\mathcal{V}}\equiv 0$,
whence $pq|_{\mathcal{V}}\equiv 0$. \ Since $pq\in \mathcal{P}_{n}$ and $%
\phi _{\beta }$ is 1-1, it follows that $(pq)(X)=0$.

(ii) \ Suppose $\beta $ is consistent. \ Let $p\in \mathcal{P}_{n}$ and let $%
q\in \mathcal{P}$, with $pq\in \mathcal{P}_{2n}$. \ If $p(X)=0$, then
clearly $\mathcal{V}_{\beta }\subseteq \mathcal{Z}(p)$, whence $(pq)|_{%
\mathcal{V}_{\beta }}\equiv 0$. \ Now, consistency implies that $\Lambda
_{\beta }(pq)=0$, so (\ref{neweq31}) holds.

Assume now that (\ref{neweq31}) holds and suppose $p,q,pq\in \mathcal{P}_{n}$
with $p(X)=0$. \ Now, for each $s\in \mathcal{P}_{n}$, $p(qs)\in \mathcal{P}%
_{2n}$, so (\ref{neweq31}) implies 
\begin{equation*}
\left\langle \mathcal{M}(n)\widehat{pq},\hat{s}\right\rangle =\Lambda
_{\beta }((pq)s)=\Lambda _{\beta }(p(qs))=0\;\;\text{(by (\ref{neweq31})).}
\end{equation*}%
Thus $(pq)(X)=\mathcal{M}(n)\widehat{pq}=0$, so $\mathcal{M}(n)$ is
recursively generated.
\end{proof}

\noindent%
It is not difficult to see that Lemma \ref{lem22} remains true if the
hypothesis that $\beta $ is consistent is replaced by the condition that $%
\phi _{\beta }$ is $1$-$1$. \ Indeed, we see that $\phi _{\beta }$ is $1$-$%
1\Leftrightarrow \ker \;\mathcal{M}(n)=\ker \;W_{n}\Leftrightarrow \operatorname{%
rank}\;\mathcal{M}(n)=\operatorname{rank}\;W_{n}$.

For the case when $\mathcal{V}\equiv \mathcal{V}_{\beta }$ is finite and the
elements of $\mathcal{V}$ can be computed exactly, we next describe an
elementary procedure for determining whether or not $\beta $ is consistent.
\ Denote the distinct points of $\mathcal{V}$ as $\{w_{j}\}_{j=1}^{m}$. \
Recall the matrix $W\equiv W_{2n}[\mathcal{V}_{\beta }]$,$\,$with $m$ rows
and with columns indexed by the monomials in $\mathcal{P}_{2n}$ (indexed, as
usual, in degree-lexicographic order). \ The entry of $W$ in row $k$, column 
$x^{i}$ is $w_{k}^{i}$. \ Clearly, a polynomial $p(x)\equiv \sum_{|i|\leq
2n}a_{i}x^{i}\in \mathcal{P}_{2n}$ satisfies $p|_{\mathcal{V}}\equiv 0$ if
and only if $W\hat{p}=0$. \ Using Gaussian elimination, we may row-reduce $W$
so as to find a basis for $\ker W$, say $\{\hat{p}_{1},...,\hat{p}_{s}\}$. \
It follows that $\{p_{1},...,p_{s}\}$ is a basis for $\{p\in \mathcal{P}%
_{2n}:p|_{\mathcal{V}}\equiv 0\}$. \ Let $\hat{p}_{j}:=(a_{ji})_{|i|\leq 2n}$
$(1\leq j\leq s)$. \ We now see that $\beta $ is consistent if, and only if,
for each $j$, $\Lambda _{\beta }(p_{j})=\sum_{|i|\leq 2n}a_{ji}\beta _{i}=0$%
. \ 

In Example \ref{ex15} (above) we were able to compute the points of $%
\mathcal{V}$ exactly and to then check the consistency of $\beta $ using the
preceding method. \ In other examples we may be able to determine that $%
\mathcal{V}$ is finite (from the form of the polynomial relations which
determine $\mathcal{V}$) without being able to exactly compute the points of
the variety. \ In such cases we cannot employ the above procedure for
checking consistency. \ 

The concluding remarks of Section \ref{sectideals}, particularly (\ref{new24}%
), suggest alternate, more algebraic, approaches to verifying consistency
that we will pursue below and in Sections \ref{sect4} and \ref{newsect6}. \
Let $\mathcal{N}_{n}:=\{p\in \mathcal{P}_{n}:\mathcal{M}(n)\hat{p}=0\}$ and
let $\mathcal{J}_{\beta }:=(\mathcal{N}_{n})$ denote the ideal of $\mathcal{P%
}$ generated by $\mathcal{N}_{n}$. \ For $\mathcal{S}\subseteq \mathcal{P}$,
let $\mathcal{V}(\mathcal{S}):=\{x\in \mathbb{R}^{d}:p(x)=0$ for every $p\in 
\mathcal{S}\}$; we have 
\begin{equation*}
\mathcal{V}_{\beta }=\mathcal{V}(\mathcal{N}_{n})=\mathcal{V}(\mathcal{J}%
_{\beta }).
\end{equation*}%
Let $\mathcal{I}(\mathcal{V}):=\{p\in \mathbb{R}^{d}[x]:p|_{\mathcal{V}%
_{\beta }}\equiv 0\}\;(\supseteq \mathcal{J}_{\beta })$, and set $\mathcal{K}%
_{n}:=\mathcal{I}(\mathcal{V})\bigcap \mathcal{P}_{n}$. \ Clearly $\mathcal{N%
}_{n}\subseteq \mathcal{K}_{n}$, and $\phi _{\beta }$ is one-to-one if and
only if $\mathcal{N}_{n}=\mathcal{K}_{n}$. \ 

Consider a polynomial ideal $\mathcal{I}\subseteq \mathcal{P}$. \ Recall
from \cite{MoSa} that $\{p_{1},...,p_{k}\}\subseteq \mathcal{I}$ forms an $H$%
\textit{-basis} for $\mathcal{I}$ if for every $p\in \mathcal{I}$ there
exist polynomials $q_{1},...,q_{k}$ such that $p=\sum_{i=1}^{k}p_{i}q_{i}$
and $\deg p_{i}q_{i}\leq \deg p\;(1\leq i\leq k)$. \ Every Gr\"{o}bner basis
is an $H$-basis; in particular, every polynomial ideal has an $H$-basis \cite%
{MoSa}. \ We will utilize the following weak $H$-basis condition for
elements of $\mathcal{P}_{2n}$:%
\begin{equation}
\begin{array}{l}
\text{For each }p\in \mathcal{I}(\mathcal{V})\text{ with }\deg p\leq 2n\text{%
, there exist }m>0\text{, } \\ 
\text{polynomials }h_{1},...,h_{m}\in \mathcal{N}_{n}\text{, and polynomials 
}f_{1},...,f_{m}\text{ } \\ 
\text{with }\deg f_{i}h_{i}\leq 2n\;(1\leq i\leq m)\text{ (where }m\text{, }%
h_{i}\text{ and }f_{i}\text{ may depend on }p\text{) } \\ 
\text{such that }p=\sum_{i=1}^{m}f_{i}h_{i}.%
\end{array}
\label{eq31new}
\end{equation}%
Note that if $\mathcal{N}_{n}$ contains an $H$-basis for $\mathcal{I}(%
\mathcal{V})$, then (\ref{eq31new}) is satisfied.

The following result is proved in \cite{Fia4}.

\begin{theorem}
\label{thm32new}(\cite{Fia4}) \ If $\mathcal{M}(n)$ is recursively generated
and satisfies (\ref{eq31new}), then $\beta ^{(2n)}$ is consistent.
\end{theorem}

\begin{corollary}
\label{cor33}If $\mathcal{M}(n)$ is recursively generated and $\mathcal{N}%
_{n}$ contains an $H$-basis for $\mathcal{I}(\mathcal{V})$, then $\beta
^{(2n)}$ is consistent.
\end{corollary}

We next present some examples which illustrate Corollary \ref{cor33}.

\begin{proposition}
\label{prop41}For $d=2$ (the plane), if $\mathcal{M}(n)(\beta )$ is
recursively generated and ${\mathcal{V}}_{\beta }$ is a proper, infinite
irreducible curve, then $\beta $ is consistent.
\end{proposition}

\begin{proof}
There is an irreducible polynomial $f\in \mathcal{P}_{n}$ such that $%
f(X,Y)=0 $ and $\mathcal{V}_{\beta }=\mathcal{Z}(f)$. \ \cite[Corollary 1,
p. 18]{Ful} implies that $\mathcal{I}(\mathcal{V})=(f)$ (the ideal generated
by $f$), and clearly $\{f\}\;(\subseteq \mathcal{N}_{n})$ is an $H$-basis
for\ $\mathcal{I}(\mathcal{V})$. \ The result now follows from Corollary \ref%
{cor33}.
\end{proof}

\begin{example}
\label{ex42}We illustrate Proposition \ref{prop41} with an example from %
\cite[Example 5.2]{tcmp9}. \ Consider the moment matrix 
\begin{equation*}
\mathcal{M}(3)=\left( 
\begin{array}{cccccccccc}
1 & 0 & 0 & 1 & 1 & 2 & 0 & 0 & 0 & 0 \\ 
0 & 1 & 1 & 0 & 0 & 0 & 3 & 1 & 1 & 2 \\ 
0 & 1 & 2 & 0 & 0 & 0 & 1 & 1 & 2 & 5 \\ 
1 & 0 & 0 & 3 & 1 & 1 & 0 & 0 & 0 & 0 \\ 
1 & 0 & 0 & 1 & 1 & 2 & 0 & 0 & 0 & 0 \\ 
2 & 0 & 0 & 1 & 2 & 5 & 0 & 0 & 0 & 3 \\ 
0 & 3 & 1 & 0 & 0 & 0 & 14 & 3 & 1 & 1 \\ 
0 & 1 & 1 & 0 & 0 & 0 & 3 & 1 & 1 & 2 \\ 
0 & 1 & 2 & 0 & 0 & 0 & 1 & 1 & 2 & 5 \\ 
0 & 2 & 5 & 0 & 0 & 3 & 1 & 2 & 5 & 33%
\end{array}%
\right) .
\end{equation*}%
It is straightforward to check that $\mathcal{M}(3)$ is positive and
recursively generated, with column relations $YX=1$, $YX^{2}=X$, $Y^{2}X=Y$
in ${\mathcal{C}}_{{\mathcal{M}}(3)}$, and $\operatorname{rank}\;\mathcal{M}(3)=7$%
. \ Then ${\mathcal{V}}_{\beta }$ is the hyperbola $yx=1$, and Proposition %
\ref{prop41} implies that $\beta $ is consistent. \ (The existence of a
representing measure for $\beta $ follows from \cite[Theorem 2.1]{tcmp9}.) \ 
\textbf{\qed}
\end{example}

Let $\mathcal{F}\equiv \{r_{1},...,r_{d}\}\subseteq \mathcal{P}\equiv 
\mathbb{R}[x_{1},...,x_{d}]$ and assume that $r_{1},...,r_{d}$ have no
common zeros at infinity. \ This means that the leading homogeneous forms $%
Lf(r_{1}),...,Lf(r_{d})$ have no common zeros except $(0,...,0)\in \mathbb{R}%
^{d}$ \cite{MoSa}. \ In this case, \cite[Theorem 5.3]{MoSa} implies that $%
\mathcal{F}$ is an $H$-basis for $\mathcal{I}:=(r_{1},...,r_{d})$, and $%
\mathcal{V}:=V_{\mathbb{R}}(\mathcal{I})$ is finite \cite[Section 7]{MoSa}.
\ Further, $H_{\mathcal{I}}(k)=\delta :=\deg r_{1}\cdot ...\cdot \deg r_{d}$
for $k\geq \deg r_{1}+...+\deg r_{d}-d+1$, and $H_{\mathcal{I}}(k)<\delta $
for $k<\deg r_{1}+...+\deg r_{d}-d+1$ \cite[Lemma 5.4]{MoSa}. \ Recall that
a common zero $w$ of $r_{1},...,r_{d}$ is simple if the Jacobian $(\frac{%
\partial r_{i}}{\partial x_{j}}(w))_{1\leq i,j\leq d}$ has rank $d$. \
Specializing to $d=2$, a theorem of M. Noether implies that if $r_{1}$ and $%
r_{2}$ have no common zeros at infinity and the common zeros are all real
and simple, then there are exactly $M:=\deg r_{1}\deg r_{2}$ common zeros, $%
\mathcal{V}=\{w_{1},...,w_{m}\}$, and if $p\in \mathcal{P}$ satisfies $p|_{%
\mathcal{V}}\equiv 0$, then $p$ has a representation $p=a_{1}r_{1}+a_{2}r_{2}
$, where $a_{i}\in \mathcal{P}$ satisfies $\deg a_{i}\leq \deg p-\deg
r_{i}\;(i=1,2)$. \ These observations, together with Corollary \ref{cor33},
lead to the following criterion for consistency. \ 

\begin{proposition}
\label{cor36}Suppose $d=2$. \ Let $\mathcal{M}(n)$ be recursively generated,
and suppose a basis for $\ker \mathcal{M}(n)$ consists of $\hat{r}_{1}$ and $%
\hat{r}_{2}$, where $r_{1}$ and $r_{2}$ have no common zeros at infinity and
whose common zeros are all real and simple. \ Then $\beta ^{(2n)}$ is
consistent. \ 
\end{proposition}

\begin{proof}
The above mentioned results of \cite{MoSa} show that $\{r_{1},r_{2}\}$ forms
and $H$-basis for $\mathcal{I}:=(r_{1},r_{2})$, and since the common zeros
of $r_{1}$ and $r_{2}$ are real and simple, Noether's Theorem implies that $%
\mathcal{I}$ coincides with $\mathcal{I}(V)$. \ The result now follows from
Corollary \ref{cor33}.
\end{proof}

\noindent%
Example \ref{ex71} (below) illustrates Proposition \ref{cor36}.

We conclude this section by illustrating a broad class of extremal moment
matrices having flat extensions (and representing measures). \ Suppose $%
\mathcal{M}(n)$ admits a positive extension $\mathcal{M}(n+1)$. \ If $f\in 
\mathcal{P}_{n}$ and $f(X)=0$ in $\mathcal{C}_{\mathcal{M}(n)}$, then $%
f(X)=0 $ in $\mathcal{C}_{\mathcal{M}(n+1)}$, i.e., $\mathcal{N}%
_{n}\subseteq \mathcal{N}_{n+1}$ \cite{FiCM}. \ If, further, $\mathcal{M}%
(n+1)$ is recursively generated, then it follows that $\mathcal{J}_{\beta
}\bigcap \mathcal{P}_{n+1}\subseteq \mathcal{N}_{n+1}$. \ Motivated by \cite%
{Moe2}, we say that $\mathcal{M}(n+1)$ is a \textit{tight} extension of $%
\mathcal{M}(n)$ if $\mathcal{N}_{n+1}=\mathcal{J}_{\beta }\bigcap \mathcal{P}%
_{n+1}$. \ (\cite{Moe2} discusses ``tight extensions'' of linear functionals
on $\mathcal{P}_{n}$.) \ 

\begin{theorem}
\label{thm36}(\cite{Fia4}) \ If $\mathcal{M}(n)\geq 0$ admits a tight flat
extension, then $\mathcal{M}(n)$ is extremal.
\end{theorem}

Recall that $\mathcal{M}(n)$ is \textit{flat} if $\operatorname{rank}\;\mathcal{M}%
(n)=\operatorname{rank}\;\mathcal{M}(n-1)$; the proof of \cite[Theorem 5.4]{tcmp1}
shows that if $\mathcal{M}(n)\;(\geq 0)$ is flat, then $\mathcal{M}(n)$
admits a tight flat extension, so $\mathcal{M}(n)$ is also extremal. \
Remarkably, examination of the proofs of \cite{tcmp6}, \cite{tcmp7}, \cite%
{tcmp9} and \cite{FiaOT} reveals that in each extremal case studied therein, 
$\mathcal{M}(n)$ admits a tight flat extension $\mathcal{M}(n+1)$. \ We can
further illustrate this phenomenon as follows.

\begin{example}
The extremal matrices $\mathcal{M}(n)$ of Example \ref{ex14} admit tight
flat extensions. \ For simplicity of notation, we consider only $\mathcal{M}%
(3)$ and $\beta \equiv \beta ^{(6)}$. \ We have $\operatorname{rank}\;\mathcal{M}%
(3)=6$, with column relations $\bar{Z}Z=1$, $\bar{Z}Z^{2}=Z$, $\bar{Z}^{2}Z=%
\bar{Z}$, and $\bar{Z}^{3}-Z^{3}=a(Z^{2}-\bar{Z}^{2})$. \ Thus $\mathcal{N}%
_{3}$ has a basis $\mathcal{B}_{3}=\{\bar{z}z-1,\bar{z}z^{2}-z,\bar{z}^{2}z-%
\bar{z},\bar{z}^{3}-z^{3}-a(z^{2}-\bar{z}^{2})\}$. \ In Example \ref{ex14}
we showed that $\gamma ^{(2n)}$ has a (unique) $\operatorname{rank}\;\mathcal{M}%
(n)$-atomic representing measure, so \cite{tcmp1}, \cite{tcmp2} imply that $%
\mathcal{M}(n)$ has a (unique, recursively generated) flat extension $%
\mathcal{M}(n+1)$. \ For the unique flat extension $\mathcal{M}(4)$ we have $%
\mathcal{N}_{4}\supseteq \mathcal{J}_{\beta }\bigcap \mathcal{P}%
_{4}\supseteq \mathcal{B}_{4}:=\mathcal{B}_{3}\bigcup $ $\{\bar{z}^{2}z^{2}-%
\bar{z}z,\bar{z}^{3}z-\bar{z}^{2},\bar{z}z^{3}-z^{2},\bar{z}%
^{3}z-z^{4}-a(z^{3}-\bar{z}^{2}z),\bar{z}^{4}-\bar{z}z^{3}-a(\bar{z}z^{2}-%
\bar{z}^{3})\}$. \ Since $\mathcal{B}_{4}$ is independent in $\mathcal{P}%
_{4} $, and $\dim \;\mathcal{N}_{4}=\dim \;\mathcal{P}_{4}-\operatorname{rank}\;%
\mathcal{M}(4)=\dim \;\mathcal{P}_{4}-\operatorname{rank}\;\mathcal{M}(3)=15-6=9$%
, we have $9=\dim \;\mathcal{N}_{4}\geq \dim \;\mathcal{J}_{\beta }\bigcap 
\mathcal{P}_{4}\geq \operatorname{card}\;\mathcal{B}_{4}=9$, whence $\mathcal{M}%
(4)$ is tight. $\ $\textbf{\qed}
\end{example}

Despite Theorem \ref{thm36} and the preceding examples, we will show in
Section \ref{newsect6} (Proposition \ref{prop61}) that there exists a
positive, recursively generated, extremal $\mathcal{M}(3)$, admitting a flat
extension, but having no tight flat extension.

\section{\label{sect3}The Extremal Moment Problem}

Assume that $\beta \equiv \beta ^{(2n)}$ is extremal, i.e., $r:=\operatorname{rank%
}\;\mathcal{M}(n)$ and $v:=\operatorname{card}\;\mathcal{V}_{\beta }$ satisfy $%
r=v $. \ Let $\mathcal{V}\equiv \{w_{1},...,w_{r}\}$ denote the distinct
points of $\mathcal{V}_{\beta }$. \ If $\mu $ is a representing measure for $%
\beta $, then $\operatorname{supp}\;\mu \subseteq \mathcal{V}$ and $r\leq 
\operatorname{card}\;\operatorname{supp}\;\mu \leq v$, so the extremal hypothesis $r=v$
implies that $\operatorname{supp}\;\mu =\mathcal{V}$. \ Thus $\mu $ is
necessarily is of the form 
\begin{equation}
\mu =\sum_{i=1}^{r}\rho _{i}\delta _{w_{i}}.  \label{eq31}
\end{equation}%
We begin by establishing a criterion which allows us to compute the
densities $\rho _{i}$.

Let $p_{1},...,p_{r}$ be polynomials in $\mathcal{P}_{n}$ such that $%
\mathcal{B}\equiv \{p_{1}(X),...,p_{r}(X)\}$ is a basis for the column space
of $\mathcal{M}(n)$, and set 
\begin{equation*}
V\equiv V_{\mathcal{B}}[\mathcal{V}]:=\left( 
\begin{array}{ccccc}
p_{1}(w_{1}) & . & . & . & p_{1}(w_{r}) \\ 
. & . & . & . & . \\ 
. & . & . & . & . \\ 
. & . & . & . & . \\ 
p_{r}(w_{1}) & . & . & . & p_{r}(w_{r})%
\end{array}%
\right) .
\end{equation*}%
Now $V$ is singular if and only if there exist scalars $\alpha
_{1},...,\alpha _{r}$, not all $0$, such that $\alpha
_{1}p_{1}(w_{i})+\cdots \alpha _{r}p_{r}(w_{i})=0$ ($1\leq i\leq r$). \
Equivalently, the polynomial $p\in \mathcal{P}_{n}$ defined by $p:=\alpha
_{1}p_{1}+\cdots +\alpha _{r}p_{r}$ satisfies $p|_{\mathcal{V}}\equiv 0$. \
Since $\mathcal{B}$ is a basis, it follows that $p(X)\equiv \alpha
_{1}p_{1}(X)+\cdots +\alpha _{r}p_{r}(X)\not=0$, so $\phi _{\beta }$ is not
1-1. \ Conversely, suppose $\phi _{\beta }$ is not 1-1, i.e., there exists $%
q\in \mathcal{P}_{n}$ with $q|_{{\mathcal{V}}}\equiv 0$ and $q(X)\not=0$ in $%
\mathcal{C}_{\mathcal{M}(n)}$. \ Since $\mathcal{B}$ is a basis, there exist
scalars $a_{1},...,a_{r}$, not all $0$, such that $q(X)=%
\sum_{i=1}^{r}a_{i}p_{i}(X)$, and since $\phi _{\beta }$ is well-defined, we
may assume that $q=\sum_{i=1}^{r}a_{i}p_{i}$. \ Now $q|_{\mathcal{V}}\equiv
0 $ implies that $\sum_{i=1}^{r}a_{i}p_{i}(w_{j})=0$ ($1\leq j\leq r$),
whence $V$ is singular. Thus we have

\begin{lemma}
\label{lem31}The following are equivalent for $\beta $ extremal: \newline
i) $\phi _{\beta }$ is 1-1, i.e., $p\in {\mathcal{P}}_{n}$, $p|_{{\mathcal{V}%
}}\equiv 0\Longrightarrow p(X)=0$ in $\mathcal{C}_{\mathcal{M}(n)}$; \newline
ii) For any basis ${\mathcal{B}}$ of $\mathcal{C}_{\mathcal{M}(n)}$, $V$ is
invertible; \newline
iii) There exists a basis ${\mathcal{B}}$ of $\mathcal{C}_{\mathcal{M}(n)}$
such that $V$ is invertible.
\end{lemma}

Suppose now that $\beta $ is extremal and let $\mathcal{B}$ be any basis for 
$\mathcal{C}_{\mathcal{M}(n)}$; thus there exist polynomials $%
p_{1},...,p_{r}\in \mathcal{P}_{n}$ such that $\mathcal{B}%
=\{p_{1}(X),...,p_{r}(X)\}$. \ If $\beta $ has a representing measure $\mu $%
, then $\phi _{\beta }$ is 1-1 (cf. Section \ref{sect1}), so Lemma \ref%
{lem31} shows that $V$ is invertible, whence $\mu $ is uniquely determined
from (\ref{eq31}) by 
\begin{equation}
(\rho _{1},...,\rho _{r})^{T}=V^{-1}(\Lambda _{\beta }(p_{1}),...,\Lambda
_{\beta }(p_{r}))^{T}.  \label{eq32}
\end{equation}%
Assuming only that $\beta $ is extremal and that $\phi _{\beta }$ is 1-1,
let $\mu _{\mathcal{B}}$ denote the measure defined by (\ref{eq31}) and (\ref%
{eq32}). \ Our main result, which follows, includes a proof of Theorem \ref%
{thm13}.

\begin{theorem}
\label{thm32}For $\beta \equiv \beta ^{(2n)}$ extremal, the following are
equivalent: \newline
(i) $\ \beta $ has a representing measure; \newline
(ii) $\ \beta $ has a unique representing measure, which is $\operatorname{rank}\;%
{\mathcal{M}}(n)$-atomic; \newline
(iii) \ For some (respectively, for every) basis ${\mathcal{B}}$ of $%
\mathcal{C}_{\mathcal{M}(n)}$, $V$ is invertible and $\mu _{{\mathcal{B}}}$
is a representing measure for $\beta $; \newline
(iv) $\ \beta $ is consistent and ${\mathcal{M}}(n)\geq 0$;\newline
(v) \ $\mathcal{M}(n)\geq 0$ has a flat extension $\mathcal{M}(n+1)$;\newline
(vi) \ $\mathcal{M}(n)\geq 0$ has a unique flat extension $\mathcal{M}(n+1)$.
\end{theorem}

\noindent%
(Note that a proof of (iv) $\Rightarrow $ (i) is contained in Remark \ref%
{rem2}; we present a different proof below.)

\begin{proof}
The implications (ii) $\Longrightarrow $ (i) $\Longrightarrow $ (iv) are
clear, so it suffices to prove (iv) $\Longrightarrow $ (iii) $%
\Longrightarrow $ (ii), and to then prove (ii) $\Longrightarrow $ (vi) $%
\Longrightarrow $ (v) $\Longrightarrow $ (i) ($\Leftrightarrow $ (ii)) . We
begin with the proof of (iv) $\Longrightarrow $ (iii). \ Let $\mathcal{B}$
be a basis for $\mathcal{C}_{\mathcal{M}(n)}$, and, as above, denote $%
\mathcal{B}\equiv \{p_{1}(X),...,p_{r}(X)\}$, where $p_{1},...,p_{r}$ are
polynomials in $\mathcal{P}_{n}$. \ Let $\mathcal{V}\equiv \mathcal{V}%
_{\beta }=\{w_{1},...,w_{r}\}$ and consider $V$ as defined above. \ Since $%
\beta $ is consistent, Proposition \ref{consistency} implies that $\phi
_{\beta }$ is 1-1, so Lemma \ref{lem31} shows that $V$ is invertible, and we
may thus consider $\mu _{\mathcal{B}}$ as defined by (\ref{eq31}) and (\ref%
{eq32}). \ To show that $\mu _{\mathcal{B}}$ is a representing measure for $%
\beta $, we first show that for $f\in \mathcal{P}_{2n}$, $\int f(x)d\mu _{%
\mathcal{B}}(x)=\Lambda _{\beta }(f)$. \ Let $v_{f}:=(f(w_{1}),...,f(w_{r}))$%
. Since $V$ is invertible, there exists $a_{f}\equiv (a_{1},...,a_{r})\in $%
\textbf{$\mathbb{R}$}$^{r}$ such that $V^{T}a_{f}^{T}=v_{f}^{T}$. \ Thus $%
p\equiv \sum_{i=1}^{r}a_{i}p_{i}\in \mathcal{P}_{n}$ satisfies $%
p(w_{i})=f(w_{i})$ $(1\leq i\leq r)$. \ Now 
\begin{align*}
\int f(x)d\mu _{\mathcal{B}}(x)& =\sum_{k=1}^{r}\rho
_{k}f(w_{k})=\sum_{k=1}^{r}\rho _{k}p(w_{k}) \\
& =\sum_{k=1}^{r}\rho
_{k}\sum_{i=1}^{r}a_{i}p_{i}(w_{k})=\sum_{i=1}^{r}a_{i}\sum_{k=1}^{r}\rho
_{k}p_{i}(w_{k}) \\
& =\sum_{i=1}^{r}a_{i}\Lambda _{\beta }(p_{i})\;\;\text{(from (\ref{eq32}))}
\\
& =\Lambda _{\beta }(\sum_{i=1}^{r}a_{i}p_{i})=\Lambda _{\beta }(f) \\
\text{(since }\beta \text{ is consistent and }f-p& \in \mathcal{P}_{2n}\text{
satisfies }(f-p)|_{\mathcal{V}}\equiv 0\text{).}
\end{align*}%
To complete the proof that $\mu _{\mathcal{B}}$ is a representing measure,
it remains to show that $\mu _{\mathcal{B}}\geq 0$. \ For $1\leq k\leq r$,
let $V_{k}\equiv V_{k}(x)$ denote the matrix obtained from $V$ by replacing $%
w_{k}$ (in column $k$) by the variable $x$, and let $f_{k}\in \mathcal{P}%
_{n} $ be defined by $f_{k}(x):=\det \;V_{k}(x)$. Clearly, $%
f_{k}(w_{j})=\delta _{kj}\det \;V$ $(1\leq k,j\leq r)$. Now 
\begin{eqnarray*}
0 &\leq &\left\langle \mathcal{M}(n)\hat{f}_{k},\hat{f}_{k}\right\rangle
=\Lambda _{\beta }(f_{k}^{2})=\int f_{k}^{2}d\mu _{\mathcal{B}}\;\;\text{%
(from the preceding paragraph)} \\
&=&\sum_{j=1}^{r}\rho _{j}f_{k}^{2}(w_{j})=\rho _{k}(\det \;V)^{2},
\end{eqnarray*}%
and since $\det \;V\not=0$, it follows that $\rho _{k}\geq 0$. (Since $%
\operatorname{card}\;\operatorname{supp}\;\mu _{\mathcal{B}}=r$, it then follows that $%
\rho _{k}>0~(1\leq k\leq r)$.)

To prove (iii) $\Longrightarrow $ (ii), assume that $\nu $ is a representing
measure for $\beta $. Since $\beta $ is extremal, $\nu $ is of the form $\nu
=\sum_{i=1}^{r}\sigma _{i}\delta _{w_{i}}$ for $\sigma _{i}>0$ $(1\leq i\leq
r)$. \ Suppose $\mathcal{B}\equiv \{p_{1}(X),...,p_{r}(X)\}$ is a basis for $%
\mathcal{C}_{\mathcal{M}(n)}$ (as above) such that $V$ is invertible and $%
\mu _{\mathcal{B}}$ is a representing measure for $\beta $. \ Since $\nu $
and $\mu _{\mathcal{B}}$ are representing measures, we have 
\begin{eqnarray*}
V(\rho _{1},...,\rho _{r})^{T} &=&(\int p_{1}d\mu _{\mathcal{B}},...,\int
p_{r}d\mu _{\mathcal{B}})^{T} \\
&=&(\Lambda _{\beta }(p_{1}),...,\Lambda _{\beta }(p_{r}))^{T}=(\int
p_{1}d\nu ,...,\int p_{r}d\nu )^{T} \\
&=&V(\sigma _{1},...,\sigma _{r})^{T},
\end{eqnarray*}%
and since $V$ is invertible, it follows that $\nu =\mu _{\mathcal{B}}$. \
This completes the equivalence of (i), (ii), (iii) and (iv).

Now recall that $\beta $ has a $\operatorname{rank}\;\mathcal{M}(n)$-atomic
representing measure if and only if $\mathcal{M}(n)\geq 0$ admits a flat
extension $\mathcal{M}(n+1)$ \cite[Theorem 5.13]{tcmp1}, and clearly
distinct flat extensions correspond to distinct $\operatorname{rank}\;\mathcal{M}%
(n)$-atomic representing measures. \ Thus we have (ii) $\Rightarrow $ (vi) $%
\Rightarrow $ (v) $\Rightarrow $ (i), and since (i) $\Leftrightarrow $ (ii),
the proof is complete.
\end{proof}

\begin{remark}
For a positive, extremal $\mathcal{M}(n)$ for which the points of the
variety are known, Theorem \ref{thm32} provides two ways to determine
whether or not $\beta $ has a representing measure. \ Following Theorem \ref%
{thm32}(iv) one can use the method of Section \ref{sect2} to determine
whether or not $\beta $ is consistent. \ Alternatively, one can select any
basis ${\mathcal{B}}$ of $\mathcal{C}_{\mathcal{M}(n)}$ and check whether $V$
is invertible. \ If $V$ is not invertible, there is no representing measure.
\ If $V$ is invertible, then $\mu _{{\mathcal{B}}}$ automatically
interpolates all moments up to degree $n$, so the proof of Theorem \ref%
{thm32}(iv) $\Rightarrow $ (iii) shows that $\beta $ has a representing
measure if and only if $\mu _{{\mathcal{B}}}$ interpolates all moments of
degrees $n+1,n+2,...,2n$, in which case $\mu _{\beta }\geq 0$. \ In a given
numerical problem, one approach or the other may be easier to implement,
depending on the size of $n$ and the value of $\operatorname{rank}\;\mathcal{M}%
(n) $.
\end{remark}

\section{\label{sect4}Solution of the $\mathcal{M}(3)$ Extremal Problem with 
$Y=X^{3}:r=v=7$}

In this section (and the next) we return to the question as to whether a
positive, extremal, recursively generated moment matrix has a representing
measure (cf., Question \ref{quest2}). \ We also consider the extent to which
recursiveness implies consistency in an extremal moment problem. \ Our
motivation is the observation that it is generally much easier to verify
recursiveness than consistency. \ We examine these issues in detail for an
extremal planar moment matrix $\mathcal{M}(3)$ with $\mathcal{M}(3)\geq 0$, $%
\mathcal{M}(2)>0$, and $Y=X^{3}$ in $\mathcal{C}_{\mathcal{M}(3)}$. \ Our
first result illustrates an extremal problem in which recursiveness does
imply consistency.

\begin{theorem}
\label{thm43}Let $d=2$. \ Suppose $Y=X^{3}$ in ${\mathcal{C}}_{{\mathcal{M}}%
(3)}$. \ If $\mathcal{M}(3)$ is positive, recursively generated, and $v=r=7$%
, then $\beta ^{(6)}$ has a unique, $7$-atomic, representing measure;
equivalently, $\beta ^{(6)}$ is consistent.
\end{theorem}

\begin{example}
\label{ex44}We illustrate Theorem \ref{thm43} with the following moment
matrix: 
\begin{equation*}
\mathcal{M}(3)=\left( 
\begin{array}{cccccccccc}
1 & 0 & 0 & 1 & 2 & 5 & 0 & 0 & 0 & 0 \\ 
0 & 1 & 2 & 0 & 0 & 0 & 2 & 5 & 14 & 42 \\ 
0 & 2 & 5 & 0 & 0 & 0 & 5 & 14 & 42 & 200 \\ 
1 & 0 & 0 & 2 & 5 & 14 & 0 & 0 & 0 & 0 \\ 
2 & 0 & 0 & 5 & 14 & 42 & 0 & 0 & 0 & 0 \\ 
5 & 0 & 0 & 14 & 42 & 200 & 0 & 0 & 0 & 0 \\ 
0 & 2 & 5 & 0 & 0 & 0 & 5 & 14 & 42 & 200 \\ 
0 & 5 & 14 & 0 & 0 & 0 & 14 & 42 & 200 & 5868 \\ 
0 & 14 & 42 & 0 & 0 & 0 & 42 & 200 & 5868 & 386568 \\ 
0 & 42 & 200 & 0 & 0 & 0 & 200 & 5868 & 386568 & 26992856%
\end{array}%
\right) .
\end{equation*}%
$\mathcal{M}(3)$ is positive and recursively generated, with column basis ${%
\mathcal{B}}:=\{1,~X,~Y,~X^{2},$ $YX,~Y^{2},~YX^{2}\}$, and column relations 
$Y=X^{3}$, $Y^{2}X=208X-282Y+74YX^{2}$, and $Y^{3}=15392X-20660Y+5194YX^{2}$%
. \ A calculation shows that ${\mathcal{V}}_{\beta }$ consists of exactly $7$
points in \textbf{$\mathbb{R}$}$^{2}$, $\{(x_{i},x_{i}^{3})\}_{i=1}^{7}$,
with $x_{1}=0$, $x_{2}\cong 8.36748$, $x_{3}\cong 0.996357$, $x_{4}\cong
1.7299$, and $x_{4+j}=-x_{j+1}$ $(1\leq j\leq 3)$. \ Thus $\beta $ is
extremal, so Theorem \ref{thm43} implies that $\beta $ has a representing
measure. \ Indeed, following the method of Section \ref{sect3}, a
calculation shows that $V_{{\mathcal{B}}}$ is invertible and that $\mu _{{%
\mathcal{B}}}$ has densities $\rho _{1}\cong 0.331731$, $\rho _{2}\cong
3.3378229\times 10^{-10}$, $\rho _{3}\cong 0.249980$, $\rho _{4}\cong
0.08415439$, and $\rho _{4+j}=\rho _{j+1}$ $(1\leq j\leq 3)$. \ \textbf{\qed}
\end{example}

We begin the proof of Theorem \ref{thm43} with some preliminary results. \
Recall from Section \ref{sect2} the map $\phi _{\beta }:\mathcal{C}_{%
\mathcal{M}(n)}\rightarrow \mathcal{P}_{n}|_{\mathcal{V}_{\beta }}$, given
by $p(X)\mapsto p|_{\mathcal{V}_{\beta }}\;(p\in \mathcal{P}_{n})$. \ As
noted in Section \ref{sect2}, $\phi _{\beta }$ is 1-1 if and only if $%
\mathcal{N}_{n}=\mathcal{K}_{n}$ (where $\mathcal{N}_{n}:=\{p\in \mathcal{P}%
_{n}:p(X)=0\}$ and $\mathcal{K}_{n}:=\mathcal{I}(\mathcal{V}_{\beta
})\bigcap \mathcal{P}_{n}=\{p\in \mathcal{P}_{n}:p|_{\mathcal{V}_{\beta
}}\equiv 0\}$); we always have $\mathcal{N}_{n}\subseteq \mathcal{K}_{n}$. \ 

\begin{lemma}
\label{lem45}If $\mathcal{M}(n)(\beta )$ satisfies $r\leq v$ and $\dim {%
\mathcal{K}}_{n}\leq \dim {\mathcal{P}}_{n}-v$, then $\mathcal{M}(n)(\beta )$
is extremal and $\phi _{\beta }$ is 1-1.
\end{lemma}

\begin{proof}
We have $v\leq \dim \mathcal{P}_{n}-\dim \mathcal{K}_{n}\leq \dim \mathcal{P}%
_{n}-\dim \mathcal{N}_{n}=r\leq v$. \ It follows that $r=v$ and $\mathcal{N}%
_{n}=\mathcal{K}_{n}$, so $\mathcal{M}(n)(\beta )$ is extremal and $\phi
_{\beta }$ is 1-1.
\end{proof}

\begin{lemma}
\label{lem46}If $\mathcal{M}(3)(\beta )$ satisfies $Y=X^{3}$ and $r\leq v=7$%
, then $\phi _{\beta }$ is 1-1.
\end{lemma}

\begin{proof}
Suppose $p(x,y)\equiv
c_{1}+c_{2}x+c_{3}y+c_{4}x^{2}+c_{5}yx+c_{6}y^{2}+c_{7}x^{3}+c_{8}yx^{2}+c_{9}y^{2}x+c_{10}y^{3} 
$ is an element of $\mathcal{K}_{3}$, i.e., $p|_{\mathcal{V}_{\beta }}\equiv
0$. \ Denote the distinct points of $\mathcal{V}_{\beta }$ by $%
\{(x_{i},y_{i})\}_{i=1}^{7}$; since $y_{i}=x_{i}^{3}\;(1\leq i\leq 7)$, the $%
x_{i}$'s are distinct. \ Consider the linear map $\Psi :\mathcal{K}%
_{3}\rightarrow \mathbb{R}^{3}$ defined by $\Psi (p)=(c_{7},c_{9},c_{10})$.
\ We claim that $\Psi $ is 1-1; for, suppose $c_{7}=c_{9}=c_{10}=0$ and let $%
f(x):=p(x,x^{3})\equiv
c_{1}+c_{2}x+c_{3}x^{3}+c_{4}x^{2}+c_{5}x^{4}+c_{6}x^{6}+c_{8}x^{5}$. \
Since $f$ has the seven distinct roots $\{x_{i}\}_{i=1}^{7}$, it follows
that $c_{1}=c_{2}=c_{3}=c_{4}=c_{5}=c_{6}=c_{8}=0$, whence $p\equiv 0$ and $%
\Psi $ is 1-1. \ Thus $\dim \mathcal{K}_{3}\leq \dim \mathbb{R}%
^{3}=3=10-7=\dim \mathcal{P}_{3}-v$, so Lemma \ref{lem45} implies that $\phi
_{\beta }$ is 1-1.
\end{proof}

\begin{proposition}
\label{prop47}Let $\mathcal{M}(3)(\beta )\geq 0$, with $Y=X^{3}$ in ${%
\mathcal{C}}_{{\mathcal{M}}(3)}$. \ If ${\mathcal{B}}:=\{1,~X,~Y,~X^{2},$ $%
YX,~Y^{2},~YX^{2}\}$ is a basis for ${\mathcal{C}}_{{\mathcal{M}}(3)}$ and $%
v=r$, then $\beta ^{(6)}$ has a representing measure.
\end{proposition}

\begin{proof}
Let $\mathcal{V}\equiv \mathcal{V}_{\beta }$; Lemmas \ref{lem31} and \ref%
{lem46} imply that $V_{\mathcal{B}}[\mathcal{V}]$ is invertible, so, as in
the proof of Theorem \ref{thm32}, to prove that $\mu _{\mathcal{B}}$ is a
representing measure, it suffices to prove that $\mu _{\mathcal{B}}$ is 
\textit{interpolating} for $\beta ^{(6)}$, i.e., $\beta _{ij}=\int
y^{j}x^{i}\;d\mu _{\mathcal{B}}$ $(i,j\geq 0,~i+j\leq 6)$. \ Relation (\ref%
{eq32}) shows that $\mu _{\mathcal{B}}$ interpolates the moments
corresponding to elements of $\mathcal{B}$, namely $\beta _{00},~\beta
_{10},~\beta _{01},~\beta _{20},~\beta _{11},~\beta _{02},$ and $\beta _{21}$%
. \ From the hypothesis, we have 
\begin{equation}
Y=X^{3}.  \label{eq41}
\end{equation}%
Also, there exist $\alpha $, $\gamma \in \mathbb{R}$ and $p,q\in \mathcal{P}%
_{2}$, such that we have column relations 
\begin{equation}
Y^{2}X=\alpha YX^{2}+p(X,Y),  \label{eq42}
\end{equation}%
and%
\begin{equation}
Y^{3}=\gamma YX^{2}+q(X,Y).  \label{eq43}
\end{equation}%
In $\operatorname{supp}\;\mu _{\mathcal{B}}$ we have $y=x^{3}$, so $\int
x^{3}d\mu _{\mathcal{B}}=\int yd\mu _{\mathcal{B}}=\beta _{01}=\left\langle
Y,1\right\rangle =\left\langle X^{3},1\right\rangle =\beta _{30}$ (by (\ref%
{eq41})); thus $\int x^{3}d\mu _{\mathcal{B}}=\beta _{30}$. \ Similarly, 
\begin{eqnarray*}
\int y^{2}x\;d\mu _{\mathcal{B}} &=&\int (\alpha yx^{2}+p(x,y))d\mu _{%
\mathcal{B}}=\alpha \beta _{21}+\Lambda _{\beta }(p) \\
&=&\left\langle \alpha YX^{2}+p(X,Y),1\right\rangle =\left\langle
Y^{2}X,1\right\rangle =\beta _{12}
\end{eqnarray*}%
(by (\ref{eq32}) and (\ref{eq42})), and 
\begin{eqnarray*}
\int y^{3}d\mu _{\mathcal{B}} &=&\int (\gamma yx^{2}+q(x,y))d\mu _{\mathcal{B%
}}=\gamma \beta _{21}+\Lambda _{\beta }(q) \\
&=&\left\langle \gamma YX^{2}+q(X,Y),1\right\rangle =\left\langle
Y^{3},1\right\rangle =\beta _{03}
\end{eqnarray*}%
(by (\ref{eq32}) and (\ref{eq43})). \ Thus, $\mu _{\mathcal{B}}$
interpolates all moments up to degree $3$.

The proof now continues inductively, using the results for all degrees $<k$
to obtain the result for degree $k$, and using (\ref{eq41})-(\ref{eq43}) in
successive rows of $\mathcal{M}(3)$. \ For example, to obtain results for
degree $4$, we start with the relations $y=x^{3}$, $y^{2}x=\alpha
yx^{2}+p(x,y)$, and $y^{3}=\gamma yx^{2}+q(x,y)$, valid in $\mathcal{V}%
_{\beta }$, to get new relations of degree $4$ in $\mathcal{V}_{\beta }$: $%
x^{4}=yx$, $yx^{3}=y^{2}$, $y^{2}x^{2}=\alpha yx^{3}+xp(x,y)$, $%
y^{3}x=\gamma yx^{3}+xq(x,y)$, $y^{4}=\gamma y^{2}x^{2}+yq(x,y)$. \ Using (%
\ref{eq41})-(\ref{eq43}) and the results for degrees $1$, $2$ and $3$, we
may now successively integrate these new relations to obtain $\beta
_{i+j}=\int y^{j}x^{i}\;d\mu _{\mathcal{B}}$ $(i,j\geq 0,~i+j=4)$; for
example, $\int x^{4}d\mu _{\mathcal{B}}=\int yxd\mu _{\mathcal{B}}=\beta
_{11}=\left\langle Y,X\right\rangle =\left\langle X^{3},X\right\rangle
=\beta _{40}$. \ Degrees $5$ and $6$ are treated similarly.
\end{proof}

\begin{proof}[Proof of Theorem \ref{thm43}]
In view of Theorem \ref{thm32}(i) $\Leftrightarrow $ (ii), it suffices to
show that $\beta ^{(6)}$ has a representing measure. \ The results in \cite%
{tcmp5}, \cite{tcmp7} and \cite{tcmp9} show that if $\mathcal{M}(n)$ is
positive, recursively generated, satisfies $r\leq v$ and has a column
relation of degree one or two, then $\beta ^{(2n)}$ admits a representing
measure. \ We may thus assume that $\mathcal{M}(2)$ is positive and
invertible; indeed, positivity in $\mathcal{M}(3)$ implies that any
dependence relation in the columns of $\mathcal{M}(2)$ extends to the
columns of $\mathcal{M}(3)$ \cite{tcmp3}. \ In particular, we may assume in
the sequel that a basis $\mathcal{B}$ of $\mathcal{C}_{\mathcal{M}(3)}$
includes $\{1,X,Y,X^{2},YX,Y^{2}\}$.

Lemma \ref{lem46} implies that $\phi _{\beta }$ is 1-1. \ As in Section \ref%
{sect3}, we may thus form $\mu _{\mathcal{B}}$, and as in the proof of
Theorem \ref{thm32}, it suffices to show that $\mu _{\mathcal{B}}$ is
interpolating for $\beta $. \ The proof of Proposition \ref{prop47} shows
that this is the case if $\mathcal{B}=\{1,~X,~Y,~X^{2},~YX,~Y^{2},~YX^{2}\}$%
. This proof shows, more generally, that $\mu _{\mathcal{B}}$ is
interpolating if $\mathcal{B}$ contains $\{1,$ $X,$ $Y,$ $X^{2},$ $YX,$ $%
Y^{2}\}$ and there exist column relations of the form (\ref{eq42}) and (\ref%
{eq43}).

We consider next the case when $\mathcal{B}%
=\{1,~X,~Y,~X^{2},~YX,~Y^{2},~Y^{2}X\}$, with column relations $%
YX^{2}=u(X,Y)+\gamma Y^{2}X~(\gamma \in \mathbb{R},~\deg \;u\leq 2)$ and $%
Y^{3}=\delta Y^{2}X+t(X,Y)\;\;(\delta \in \mathbb{R},\deg \;t\leq 2)$. \ Let 
$h(x,y):=x^{2}y-u(x,y)-\gamma xy^{2}$, so that $h(X,Y)=0$ and $\mathcal{V}%
_{\beta }\subseteq \{(x,y)\in \mathbb{R}^{2}:y=x^{3}$ and $h(x,y)=0\}$. \ If 
$\gamma =0$, then $h(x,y)=0$ has at most $6$ real roots of the form $%
(x,x^{3})$, contradicting $r=v=7$. Thus $\gamma \not=0$, and we may derive a
system as in (\ref{eq42})-(\ref{eq43}); indeed, $Y^{3}=\frac{\delta }{\gamma 
}YX^{2}+(t-\frac{\delta }{\gamma }u)(X,Y)$. \ Using this system, we may now
proceed as in the proof of Proposition \ref{prop47} to conclude that $\mu _{%
\mathcal{B}}$ is interpolating. \ Finally, we consider the case $\mathcal{B}%
=\{1,~X,~Y,~X^{2},~YX,~Y^{2},~Y^{3}\}$, with relations 
\begin{equation}
YX^{2}=s(X,Y)+\delta Y^{3}\;\;(\delta \in \mathbb{R},~\deg s\leq 2)
\label{eq46}
\end{equation}%
and 
\begin{equation}
Y^{2}X=t(X,Y)+\epsilon Y^{3}\;\;(\epsilon \in \mathbb{R},~\deg t\leq 2).
\label{eq47}
\end{equation}%
Since $h(x,y):=yx^{2}-s(x,y)$ has at most $6$ roots of the form $(x,x^{3})$,
then $v=7$ implies $\delta \not=0$. \ We may now successively transform (\ref%
{eq46}) and (\ref{eq47}) into (\ref{eq42}) and (\ref{eq43}) and then apply
the method of the proof of Proposition \ref{prop47}.
\end{proof}

\section{\label{newsect6}The Extremal Problem for $\mathcal{M}(3)$ with $%
Y=X^{3}:r=v=8$}

In this section we study the extremal moment problem for a moment matrix $%
\mathcal{M}(3)$ satisfying 
\begin{equation}
\mathcal{M}(3)\geq 0,\mathcal{M}(2)>0,Y=X^{3}\text{ in }\mathcal{C}_{%
\mathcal{M}(3)}\text{, and }r=v=8.  \label{eq61new}
\end{equation}%
In Proposition \ref{prop61} we illustrate (\ref{eq61new}) with the first
example of an extremal moment matrix $\mathcal{M}(n)$, which admits a
representing measure, but for which (i) the ideal $\mathcal{J}_{\beta }$
corresponding to $\ker \;\mathcal{M}(n)$ is \textit{not} a real ideal, and
(ii) the unique flat extension $\mathcal{M}(n+1)$ is \textit{not} a tight
flat extension. \ In Theorem \ref{thm62} we resolve Question \ref{quest2} in
the negative, by constructing a moment matrix $\mathcal{M}(3)$ which
satisfies (\ref{eq61new}), but is not consistent, and thus admits no
representing measure. \ In Theorem \ref{thm63} we provide a simplified
consistency test for moment matrices satisfying (\ref{eq61new}), thereby
completing the analysis of the extremal moment problem for $\mathcal{M}(3)$
with $Y=X^{3}$ (cf. Remark \ref{rem65}(iii)). \ 

We begin by introducing the objects that we will use in our examples. \ Let $%
f(x,y):=y-x^{3}$. \ Recall from Bezout's Theorem (\cite[Theorem 8.7.10]%
{CLO'S} that if $\deg \;g=3$, then $f$ and $g$ have exactly $9$
common zeros (counting multiplicity), including complex zeros and zeros at
infinity. \ To construct a variety that will serve as $\mathcal{V}(\mathcal{M%
}(3))$ in Proposition \ref{prop61} and Theorem \ref{thm62}, we first seek a
polynomial $g\in \mathbb{R}[x,y]$ of degree $3$ such that $f$ and $g$ have
exactly $8$ distinct common real affine zeros, one of which is a zero of
multiplicity $2$. \ For this, let $\ell _{i}(x,y)=0\;(i=1,2,3)$ be lines in
the plane such that $\ell _{1}$ intersects $y=x^{3}$ in $3$ distinct points (%
$(x_{i},y_{i}),1\leq i\leq 3$), $\ell _{2}$ intersects $y=x^{3}$ in $3$
additional distinct points ($(x_{i},y_{i}),4\leq i\leq 6$), and $\ell _{3}$
intersects $y=x^{3}$ in $2$ additional distinct points ($(x_{i},y_{i}),7\leq
i\leq 8$), such that $\ell _{3}$ is the tangent line to $y=x^{3}$ at $%
(x_{8},y_{8})$. \ Setting $g(x,y):=\ell _{1}(x,y)\ell _{2}(x,y)\ell
_{3}(x,y) $, we have $\mathcal{V}((f,g))=\{(x_{i},y_{i})\}_{i=1}^{8}$, and $%
(x_{8},y_{8})$ is a common zero of $f$ and $g$ with multiplicity $2$. \
Indeed, $(x_{8},y_{8})$ is a multiple zero since $\ell _{3}(x,y)=0$ is a
common tangent line for $f(x,y)=0$ and $g(x,y)=0$ at $(x_{8},y_{8})$;
equivalently, there exist $a,b\in \mathbb{R}$ such that the differential $D:%
\mathcal{P}\rightarrow \mathbb{R}$ defined by 
\begin{equation}
D(p):=a\frac{\partial p}{\partial x}(x_{8},y_{8})+b\frac{\partial p}{%
\partial y}(x_{8},y_{8})  \label{Deq}
\end{equation}%
satisfies $D(f)=D(g)=0$ (cf. \cite[Proposition 3.4.2]{CLO'S}, \cite{MMM}). \
We next introduce some ideals which will be referenced in the sequel. \ Let $%
\mathcal{V}\equiv \mathcal{V}((f,g))\;(=\{(x_{i},y_{i})\}_{i=1}^{8})$ and
set $\mathcal{A}:=\mathcal{I}(\mathcal{V})\equiv \{p\in \mathcal{P}:p|_{%
\mathcal{V}}\equiv 0\}$ and $\mathcal{D}:=\{p\in \mathcal{A}:D(p)=0\}$; $%
\mathcal{A}$ is a real ideal (cf. Section \ref{sectideals}), and $\mathcal{D}
$ is an ideal (which contains $f$ and $g$). \ For the last assertion, note
that if $p\in \mathcal{D}$ and $q\in \mathcal{P}$, then $(pq)|_{\mathcal{V}%
}\equiv 0$ and $D(pq)=q(x_{8},y_{8})D(p)+p(x_{8},y_{8})D(q)=0$ (since $%
D(p)=0 $ and $p|_{\mathcal{V}}\equiv 0$).

As we show below, the conditions of (\ref{eq61new}) imply that $\mathcal{B}%
:=\{1,X,Y,X^{2},YX,Y^{2},$ $YX^{2},Y^{2}X\}$ is a basis for $\mathcal{C}_{%
\mathcal{M}(3)}$, so we will further require that the points of $\mathcal{V}$
are in ``general position'' relative to the monomials $%
1,x,y,x^{2},yx,y^{2},yx^{2}$ and $y^{2}x$, i.e., we will require that $%
V\equiv V_{\mathcal{B}}[\mathcal{V}]$ is invertible (cf. Lemma \ref{lem31}).
\ Let $W\equiv W_{\mathcal{B}}[\mathcal{V}]:=V^{T}$. \ Now, if $H(x,y)$ is
any real-valued function defined on $\mathcal{V}$, then there exist scalars $%
\alpha _{1},...,\alpha _{8}\in \mathbb{R}$ such that 
\begin{equation*}
H(x,y)=\alpha _{1}+\alpha _{2}x+\alpha _{3}y+\alpha _{4}x^{2}+\alpha
_{5}yx+\alpha _{6}y^{2}+\alpha _{7}yx^{2}+\alpha _{8}y^{2}x\;\;((x,y)\in 
\mathcal{V});
\end{equation*}%
indeed, $\alpha \equiv (\alpha _{1},...,\alpha _{8})$ is uniquely determined
from 
\begin{equation}
\alpha ^{T}=W^{-1}(H(x_{1},y_{1}),..,H(x_{8},y_{8}))^{T}.  \label{eq62new}
\end{equation}%
In particular, there exist unique real numbers $a_{1},...,a_{8}$ such that 
\begin{equation}
h(x,y):=y^{2}x^{2}-(a_{1}+a_{2}x+a_{3}y+a_{4}x^{2}+a_{5}yx+a_{6}y^{2}+a_{7}yx^{2}+a_{8}y^{2}x)
\label{eq62prime}
\end{equation}%
vanishes on $\mathcal{V}$.

For the sake of definiteness, let 
\begin{align}
& \ell _{1}(x,y)\text{ }\mathbf{%
:=%
}\text{ }y-4x  \notag \\
((x_{1},y_{1})& =(-2,-8),(x_{2},y_{2})=(0,0),(x_{3},y_{3})=(2,8)),  \notag \\
&  \notag \\
& \ell _{2}(x,y)\text{ }\mathbf{%
:=%
}\text{ }y-4x+3  \notag \\
((x_{4},y_{4})& =(1,1),(x_{5},y_{5})=(-\frac{1}{2}+\frac{\sqrt{13}}{2},-5+2%
\sqrt{13}),  \label{eq65new10} \\
(x_{6},y_{6})& =(-\frac{1}{2}-\frac{\sqrt{13}}{2},-5-2\sqrt{13})),  \notag \\
&  \notag \\
& \ell _{3}(x,y)\text{ }\mathbf{%
:=%
}\text{ }y-\frac{3}{4}x+\frac{1}{4}  \notag \\
((x_{7},y_{7})& =(-1,-1),(x_{8},y_{8})=(\frac{1}{2},\frac{1}{8})).  \notag
\end{align}%
Then 
\begin{eqnarray}
\;\;\;\;\;\;\;\;\;\;4g(x,y) &=&4(y-4x)(y-4x+3)(y-\frac{3}{4}x+\frac{1}{4})
\label{eq63newb} \\
&=&-48x^{3}+88yx^{2}-35y^{2}x+4y^{3}+52x^{2}-65yx+13y^{2}-12x+3y.  \notag
\end{eqnarray}%
A calculation shows that $\ell _{3}$ is tangent to both $f$ and $g$ at $%
(x_{8}y_{8})$; indeed, $D(f)=D(g)=0$, where $D$ is the functional given by (%
\ref{Deq}) with $a=1,b=\frac{3}{4}$. \ Further, $\det V=\frac{98415}{4}\sqrt{%
13}\;(\neq 0)$, so $\operatorname{rank}\;V=8$. \ Applying (\ref{eq62new}) with $%
H(x,y)=y^{2}x^{2}$, we see that in (\ref{eq62prime}) we have 
\begin{equation}
h(x,y)=y^{2}x^{2}+6x-14x^{2}-\frac{11}{2}y+\frac{43}{2}yx-yx^{2}-\frac{17}{2}%
y^{2}+\frac{1}{2}y^{2}x,  \label{eq64new}
\end{equation}%
and $h|_{\mathcal{V}}\equiv 0$. \ A calculation shows that $D(h)=-\frac{405}{%
128}\;(\neq 0)$, so $h\in (\mathcal{A}\bigcap \mathcal{P}_{4})\;\backslash \;%
\mathcal{D}$.

\begin{proposition}
\label{prop61}Let $\mu :=\sum_{i=1}^{8}\delta _{(x_{i},y_{i})}$ (with $%
(x_{i},y_{i})$ from (\ref{eq65new10})) and let $\mathcal{M}(3):=\mathcal{M}%
(3)[\mu ]$. \ Then $\mathcal{M}(3)$ satisfies (\ref{eq61new}) and has the
following additional properties:\newline
(i) \ The ideal $\mathcal{J}_{\beta ^{(6)}}$ generated by $\mathcal{N}%
_{3}\equiv \{p\in \mathcal{P}_{3}:\mathcal{M}(3)\hat{p}=0\}$ is not a real
ideal;\newline
(ii) \ $\mathcal{M}(3)$ has a flat extension, but $\mathcal{M}(3)$ does not
admit a tight flat extension.
\end{proposition}

\begin{proof}
A direct calculation using the points in (\ref{eq65new10}) shows that $V_{%
\mathcal{B}}[\mathcal{V}]$ is invertible, so it follows as in the proof of
Lemma \ref{lem31} that $\mathcal{B}$ is independent in $\mathcal{C}_{%
\mathcal{M}(3)}$. \ Also, since $\mu $ is a representing measure for $%
\mathcal{M}(3)[\mu ]$, $\operatorname{supp}\;\mu =\mathcal{V}$, $f|_{\mathcal{V}%
}\equiv 0$ and $g|_{\mathcal{V}}\equiv 0$, we have $Y=X^{3}$ and $g(X,Y)=0$
in $\mathcal{C}_{\mathcal{M}(3)}$, whence $\mathcal{B}$ is a basis for $%
\mathcal{C}_{\mathcal{M}(3)}$, $\operatorname{rank}\;\mathcal{M}(3)=8$, and $%
\mathcal{V}(\mathcal{M}(3))=\mathcal{Z}(f)\bigcap \mathcal{Z}(g)=\mathcal{V}$%
. \ Thus, $\mathcal{M}(3)$ satisfies (\ref{eq61new}). \ 

(i) \ Let $\mathcal{J}\equiv \mathcal{J}_{\beta ^{(6)}}$ denote the ideal
generated by $\{p\in \mathcal{P}_{3}:\mathcal{M}(3)\hat{p}=0\}$, so that $%
\mathcal{J}=(f,g)$. \ We claim that $\mathcal{J}$ is not a real ideal. \
For, otherwise, there would exist $G\subseteq \mathbb{R}^{2}$ such that for $%
p\in \mathcal{P}$, $p|_{G}\equiv 0\iff p\in \mathcal{J}$ (cf. Section \ref%
{sectideals}). \ In this case, since $f^{2}+g^{2}\in \mathcal{J}$, then $%
(f^{2}+g^{2})|_{G}\equiv 0$, whence $G\subseteq \mathcal{V}$. \ Recall that
the function $h$ given by (\ref{eq64new}) satisfies $h|_{\mathcal{V}}\equiv
0 $ and $D(h)\neq 0$. \ Since $p\in \mathcal{J}\Rightarrow D(p)=0$, we see
that $h\notin \mathcal{J}$; but since $h|_{\mathcal{V}}\equiv 0$, then $%
h|_{G}\equiv 0$, contradicting the defining property of $G$. \ Thus, $%
\mathcal{J}$ is not a real ideal.

(ii) \ Since $\mathcal{M}(3)$ is extremal and has a representing measure
(that is, $\mu $), it has a unique flat extension $\mathcal{M}(4)$, namely $%
\mathcal{M}(4)[\mu ]$. \ Since $h|_{\mathcal{V}}\equiv 0$, we have $h(X,Y)=0$
in $\mathcal{C}_{\mathcal{M}(4)}$ \cite[Proposition 3.1]{tcmp1}, so $h\in 
\mathcal{N}_{4}$. \ Since we have shown in the proof of (i) that $h\notin 
\mathcal{J}$, we must have $\mathcal{J}\bigcap \mathcal{P}_{4}\neq \mathcal{N%
}_{4}$, so $\mathcal{M}(4)$ is not a tight flat extension.
\end{proof}

We next present an example of $\mathcal{M}(3)$ satisfying (\ref{eq61new}),
but not consistent, so that $\beta ^{(6)}$ has no representing measure; this
provides a negative answer to Question \ref{quest2}. \ We define a linear
functional $L:\mathcal{P}_{6}\rightarrow \mathbb{R}$ by 
\begin{equation}
L(p):=a_{0}D(p)+\sum_{i=1}^{8}a_{i}p(x_{i},y_{i})\;\;(p\in \mathcal{P}_{6})
\label{eq65new}
\end{equation}%
(with $D$ and $\{(x_{i},y_{i})\}_{i=1}^{8}$ as defined just previous to
Proposition \ref{prop61}, and $a_{i}\in \mathbb{R}$ $(0\leq i\leq 8)$). \
Let $\beta ^{(6)}$ be the sequence corresponding to $L$, i.e., $\beta
_{ij}:=L(x^{i}y^{j})\;(i,j\geq 0,i+j\leq 6)$. \ Let $M\equiv \mathcal{M}(3)$
be the corresponding moment matrix, which is real symmetric since 
\begin{equation*}
\left\langle M\widehat{x^{i}y^{j}},\widehat{x^{k}y^{\ell }}\right\rangle
=L(x^{i+k}y^{j+\ell })=\left\langle M\widehat{x^{k}y^{\ell }},\widehat{%
x^{i}y^{j}}\right\rangle .
\end{equation*}%
Recall $f(x,y):=y-x^{3}$ and note that $f(X,Y)=0$ in $\mathcal{C}_{\mathcal{M%
}(3)}$. \ Indeed, for $p\in \mathcal{P}_{3}$, 
\begin{equation*}
\left\langle f(X,Y),\hat{p}\right\rangle =\left\langle \mathcal{M}(3)\hat{f},%
\hat{p}\right\rangle =L(fp)=0
\end{equation*}%
(since $D(f)=0$ and $f|_{\mathcal{V}}\equiv 0$). \ Similarly, for $g$ as
defined earlier, since $D(g)=0$ and $g|_{\mathcal{V}}\equiv 0$, we have $%
g(X,Y)=0$. \ For the sake of definiteness, let $a_{i}:=1\;(0\leq i\leq 7)$.

\begin{theorem}
\label{thm62}There exists $\alpha \;(\cong 6.97093)$ such that if $%
a_{8}>\alpha $, then $\mathcal{M}(3)$ satisfies (\ref{eq61new}) (and is thus
positive, recursively generated, and extremal), but $\beta ^{(6)}$ has no
representing measure. \ In particular, $\phi _{\beta }$ is $1$-$1$, but $%
\beta $ is not consistent.
\end{theorem}

\begin{proof}
Consider $\mathcal{B}:=\{1,X,Y,X^{2},YX,Y^{2},YX^{2},Y^{2}X\}$. \ Since $%
Y=X^{3}$ and $g(X,Y)$ $=0$, $\mathcal{B}$ spans $\mathcal{C}_{\mathcal{M}%
(3)} $. \ It follows from Smul'jan's Theorem that $\mathcal{M}(3)$ is
positive semi-definite if and only if $M_{\mathcal{B}}$, the compression of $%
\mathcal{M}(3)$ to rows and columns indexed by $\mathcal{B}$, is positive
semi-definite. \ Calculating nested determinants, we see that $M_{\mathcal{B}%
}$ is positive definite if and only if $a_{8}>\alpha $, where $\alpha :=%
\frac{6012817451}{862617600}$. \ In this case, since $M_{8}>0$ and $%
f(X,Y)=0=g(X,Y)$, it follows that $\operatorname{rank}\;\mathcal{M}(3)=8$ and $%
\mathcal{V}(\mathcal{M}(3))=\mathcal{Z}(f)\bigcap \mathcal{Z}(g)=\mathcal{V}$%
. \ In particular, $\mathcal{M}(3)$ satisfies (\ref{eq61new}) (and is thus
also recursively generated). \ Further, $\phi _{\beta }$ is $1$-$1$ (see the
proof of Proposition \ref{prop61}, or use Lemma \ref{lem64} below). \ We
claim that $\beta ^{(6)}$ is not consistent. \ Indeed, the Riesz functional
for $\beta ^{(6)}$ is $L$. \ The function $h$ from (\ref{eq63newb})
satisfies $h|_{\mathcal{V}}\equiv 0$ and $D(h)\neq 0$, whence $L(h)=D(h)\neq
0$. \ Now $\beta $ is not consistent and thus has no representing measure.
\end{proof}

In view of Theorem \ref{thm32}, the existence of a representing measure in
the extremal moment problem (\ref{eq61new}) is equivalent to establishing
that the Riesz functional $\Lambda _{\beta }$ vanishes on a basis for $%
\mathcal{P}_{6}\bigcap \mathcal{I}(\mathcal{V})$, and we will show below
that $\dim \;\mathcal{P}_{6}\bigcap \mathcal{I}(\mathcal{V})=20$. \ The
substance of the next result is that, following (\ref{new24}) and the
remarks following (\ref{new24}), the test for consistency in (\ref{eq61new})
can be reduced to checking that $\Lambda _{\beta }(h)=0$ for $h$ given by (%
\ref{eq62prime}). \ 

\begin{theorem}
\label{thm63}Suppose $\mathcal{M}(3)$ satisfies (\ref{eq61new}), with $%
\mathcal{V}(\mathcal{M}(3))=\{(x_{i},y_{i})\}_{i=1}^{8}$ and column basis $%
\mathcal{B}:=\{1,X,Y,X^{2},YX,Y^{2},YX^{2},Y^{2}X\}$. \ Let $h$ be as in (%
\ref{eq62prime}). \ Then $\beta ^{(6)}$ has a representing measure if and
only if $\Lambda _{\beta }(h)=0$.
\end{theorem}

We require the following preliminary result. \ 

\begin{lemma}
\label{lem64}If $Y=X^{3}$ in $\mathcal{C}_{\mathcal{M}(3)}$ and $r\leq v=8$,
then $\phi _{\beta }$ is 1-1.
\end{lemma}

\begin{proof}
Let $\mathcal{V}\equiv \mathcal{V}(\mathcal{M}(3))$ and for $f\in \mathcal{K}%
_{3}:=\mathcal{P}_{3}\bigcap \mathcal{I}(\mathcal{V})$, write 
\begin{eqnarray}
f(x,y) &\equiv &a_{1}+a_{2}x+a_{3}y+a_{4}x^{2}+a_{5}yx+a_{6}y^{2}
\label{eq66new} \\
&&+a_{7}x^{3}+a_{8}yx^{2}+a_{9}y^{2}x+a_{10}y^{3}.  \notag
\end{eqnarray}%
Define a linear map $\Psi :\mathcal{K}_{3}\rightarrow \mathbb{R}^{2}$ by $%
\Psi (f):=(a_{7},a_{10})$. \ We claim that $\Psi $ is 1-1. \ Suppose $%
a_{7}=a_{10}=0$ and define 
\begin{equation*}
p(x):=f(x,x^{3})\equiv
a_{1}+a_{2}x+a_{4}x^{2}+a_{3}x^{3}+a_{5}x^{4}+a_{8}x^{5}+a_{6}x^{6}+a_{9}x^{7}.
\end{equation*}%
Since $\mathcal{V}\subseteq \mathcal{Z}(y-x^{3})$, the eight points of $%
\mathcal{V}$ have distinct $x$-coordinates, and $f|_{\mathcal{V}}\equiv 0$,
it follows that $p$ has at least $8$ distinct real roots. \ Since $\deg
\;p\leq 7$, we must have $a_{1}=a_{2}=a_{3}=a_{4}=a_{5}=a_{6}=a_{8}=a_{9}=0$%
, whence $f\equiv 0$, so $\Psi $ is 1-1. \ Now $\dim \mathcal{K}_{3}\leq
\dim \mathbb{R}^{2}=10-8=\dim \mathcal{P}_{3}-v$, so Lemma \ref{lem45}
implies that $\phi _{\beta }$ is 1-1.
\end{proof}

\begin{proof}[Proof of Theorem \ref{thm63}]
If $\beta \equiv \beta ^{(6)}$ has a representing measure, then $\beta $ is
consistent, and since $h\in \mathcal{P}_{6}\bigcap \mathcal{I}(\mathcal{V})$%
, it follows that $\Lambda _{\beta }(h)=0$. \ For the converse, we suppose
that $\Lambda _{\beta }(h)=0$ and we will show that $\beta $ is consistent,
i.e., $\mathcal{P}_{6}\bigcap \mathcal{I}(\mathcal{V})\subseteq \ker
\;\Lambda _{\beta }$ (cf. Theorem \ref{thm32}). \ To this end, we first
compute $\dim (\mathcal{P}_{6}\bigcap \mathcal{I}(\mathcal{V}))$. \ Consider 
$W\equiv W_{6}[\mathcal{V}]$ (cf. Section \ref{sectideals}); clearly, $p\in 
\mathcal{P}_{6}\bigcap \mathcal{I}(\mathcal{V})\iff \hat{p}\in \ker \;W$, so 
$\dim (\mathcal{P}_{6}\bigcap \mathcal{I}(\mathcal{V}))=\dim \;\ker \;W=\dim
\;\mathcal{P}_{6}-\operatorname{rank}\;W$. \ Lemma \ref{lem64} shows that $\phi
_{\beta }$ is 1-1, so Lemma \ref{lem31} implies that $W_{\mathcal{B}}[%
\mathcal{V}]\;(\equiv V_{\mathcal{B}}[\mathcal{V}]^{T})$ is invertible. \
Now $W_{\mathcal{B}}[\mathcal{V}]$ is the compression of $W$ to columns
indexed by the monomials corresponding to elements of $\mathcal{B}$, so $%
8\geq $ row\ rank $W=\operatorname{rank}\;W\geq \operatorname{rank}\;W_{\mathcal{B}}[%
\mathcal{V}]=8$, whence $\operatorname{rank}\;W=8$. \ Thus, 
\begin{equation*}
\dim (\mathcal{P}_{6}\bigcap \mathcal{I}(\mathcal{V}))=\dim \;\mathcal{P}%
_{6}-\operatorname{rank}\;W=28-8=20.
\end{equation*}

Let $f(x,y):=y-x^{3}$, so that $\mathcal{M}(3)\hat{f}=f(X,Y)=0$ in $\mathcal{%
C}_{\mathcal{M}(3)}$. \ Also, there exist $a,b\in \mathbb{R}$ and $p\in 
\mathcal{P}_{2}$ such that $g(x,y):=y^{3}+ayx^{2}+by^{2}x+p(x,y)$ satisfies $%
\mathcal{M}(3)\hat{g}=g(X,Y)=0$ in $\mathcal{C}_{\mathcal{M}(3)}$. \ Since $%
r=8$, it follows that $\mathcal{V}=\mathcal{Z}(f)\bigcap \mathcal{Z}(g)$,
and clearly $f,g\in \mathcal{I}(\mathcal{V})$. \ Now, if $s,t\in \mathcal{P}%
_{3}$, then $sf+tg\in \mathcal{P}_{6}\bigcap \mathcal{I}(\mathcal{V})$ and $%
\Lambda _{\beta }(sf+tg)=\left\langle \mathcal{M}(3)\hat{f},\hat{s}%
\right\rangle +\left\langle \mathcal{M}(3)\hat{g},\hat{t}\right\rangle =0$,
whence $sf+tg\in \ker \;\Lambda _{\beta }$ (see the remarks following (\ref%
{new24})).

We next identify $19$ linearly independent elements of $\mathcal{P}%
_{6}\bigcap \mathcal{I}(\mathcal{V})$ of the form $sf+tg\;(s,t\in \mathcal{P}%
_{3})$. \ Consider the following $20$ polynomials:

\begin{equation*}
\begin{array}{cc}
f_{1}:=g\equiv y^{3}+ayx^{2}+by^{2}x+p(x,y) & f_{2}:=f\equiv y-x^{3} \\ 
f_{3}:=xg\equiv y^{3}x+ayx^{3}+by^{2}x^{2}+xp(x,y) & f_{4}:=xf\equiv yx-x^{4}
\\ 
f_{5}:=yg\equiv y^{4}+ay^{2}x^{2}+by^{3}x+yp(x,y) & f_{6}:=yf\equiv
y^{2}-yx^{3} \\ 
f_{7}:=x^{2}g\equiv y^{3}x^{2}+ayx^{4}+by^{2}x^{3}+x^{2}p(x,y) & 
f_{8}:=x^{2}f\equiv yx^{2}-x^{5} \\ 
f_{9}:=yxg\equiv y^{4}x+ay^{2}x^{3}+by^{3}x^{2}+yxp(x,y) & f_{10}:=yxf\equiv
y^{2}x-yx^{4} \\ 
f_{11}:=y^{2}g\equiv y^{5}+ay^{3}x^{2}+by^{4}x+y^{2}p(x,y) & 
f_{12}:=y^{2}f\equiv y^{3}-y^{2}x^{3} \\ 
f_{13}:=x^{3}g\equiv y^{3}x^{3}+ayx^{5}+by^{2}x^{4}+x^{3}p(x,y) & 
f_{14}:=x^{3}f\equiv yx^{3}-x^{6} \\ 
f_{15}:=yx^{2}g\equiv y^{4}x^{2}+ay^{2}x^{4}+by^{3}x^{3}+yx^{2}p(x,y) & 
f_{16}:=yx^{2}f\equiv y^{2}x^{2}-yx^{5} \\ 
f_{17}:=y^{2}xg\equiv y^{5}x+ay^{3}x^{3}+by^{4}x^{2}+y^{2}xp(x,y) & 
f_{18}:=y^{2}xf\equiv y^{3}x-y^{2}x^{4} \\ 
f_{19}:=y^{3}g\equiv y^{6}+ay^{4}x^{2}+by^{5}x+y^{3}p(x,y) & 
f_{20}:=y^{3}f\equiv y^{4}-y^{3}x^{3}.%
\end{array}%
\end{equation*}%
We assert that $\mathcal{F}:=\{f_{i}\}_{i=1}^{19}$ is linearly independent
in $\mathcal{P}_{6}\bigcap \mathcal{I}(\mathcal{V})\bigcap \ker \;\Lambda
_{\beta }$. \ For $1\leq k\leq 19$, set $\mathcal{F}_{k}:=\{f_{i}%
\}_{i=1}^{k} $. \ Proceeding inductively, let $2\leq k\leq 19$ and assume
that $\mathcal{F}_{k-1}$ is linearly independent. \ Observe that, except
when $k=6,10,12,16,18$, $f_{k}$ contains a monomial of highest degree that
does not appear in any polynomial in $\mathcal{F}_{k-1}$, whence $\mathcal{F}%
_{k}$ is linearly independent. \ In the remaining cases, note that \newline
(i) $\ f_{k}$ contains a monomial of highest degree that also appears, among
the elements of\ $\mathcal{F}_{k}$, only in $f_{k-2}$;\newline
(ii) \ $f_{k-2}$ has a different monomial that also appears, among the
elements of $\mathcal{F}_{k}$, only in $f_{k-1}$;\newline
(iii) \ $f_{k-1}$ has a monomial of highest degree that appears in no other
element of $\mathcal{F}_{k}$.\newline
We thus see that $\mathcal{F}_{k}$ is independent in these cases. \ (Observe
also that $\mathcal{F}_{20}$ is dependent, since 
\begin{equation*}
f_{20}=-f_{13}-af_{16}-bf_{18}+f_{5}+p(y-x^{3}),
\end{equation*}%
and $p(y-x^{3})=pf_{2}\in \left\langle
f_{2},f_{4},f_{6},f_{8},f_{10},f_{12}\right\rangle $.)

Now,\ $\dim [\mathcal{P}_{6}\bigcap \mathcal{I}(\mathcal{V})\bigcap \ker
\;\Lambda _{\beta }]\geq 19$ and $\dim (\mathcal{P}_{6}\bigcap \mathcal{I}(%
\mathcal{V}))=20$. \ Since $h\in \mathcal{P}_{6}\bigcap \mathcal{I}(\mathcal{%
V})\bigcap \ker \;\Lambda _{\beta }$, to complete the proof that $\mathcal{P}%
_{6}\bigcap \mathcal{I}(\mathcal{V})\subseteq \ker \;\Lambda _{\beta }$, it
suffices to verify that $h\notin \left\langle
\{f_{i}\}_{i=1}^{19}\right\rangle $. \ Let $2\leq k\leq 19$ and assume by
induction that $h\notin \left\langle \{f_{i}\}_{i=1}^{k-1}\right\rangle $. \
Consider a linear combination $q:=\alpha _{1}f_{1}+\cdots +\alpha _{k}f_{k}$%
, with $\alpha _{k}\neq 0$. \ Except when $k=6,10,12,16,18$, $f_{k}$
contains a monomial term of highest degree that does not appear in $h$ or in
any element of $\mathcal{F}_{k-1}$, so $q\neq h$. \ In the remaining cases,
if $q=h$, then proceeding as in the proof that $\mathcal{F}$ is independent,
we see that $\alpha _{k-2}\neq 0$, and then that $\alpha _{k-1}\neq 0$. \
Now $f_{k-1}$ contains a monomial of highest degree that does not appear in $%
h$ or in any other element of $\mathcal{F}_{k}$, so we arrive at a
contradiction. \ Thus $q\neq h$ in these cases also. \ Now, following (\ref%
{new24}), $\{f_{i}\}_{i=1}^{19}\bigcup \{h\}$ forms a basis for $\mathcal{P}%
_{6}\bigcap \mathcal{I}(\mathcal{V})\bigcap \ker \;\Lambda _{\beta }$,
whence $\mathcal{P}_{6}\bigcap \mathcal{I}(\mathcal{V})\subseteq \ker
\;\Lambda _{\beta }$, so $\beta $ is consistent. \ The proof is now complete.
\end{proof}

\begin{remark}
\label{rem65}(i) \ If the points of the variety $\{(x_{i},y_{i})\}_{i=1}^{8}$
are known explicitly, then $h$ can be computed as in (\ref{eq62prime}). \ In
this case, Theorem \ref{thm63} provides an effective test for the existence
of a representing measure in the extremal problem (\ref{eq61new}). \ If the
points of the variety are not known explicitly, there is still available a
concrete test for the non-existence of a representing measure, as follows:%
\newline
If $\mathcal{M}(3)$ (as in (\ref{eq61new})) has a representing measure, then
there is a unique flat extension $\mathcal{M}(4)$, and $\mathcal{V}(\mathcal{%
M}(4))=\mathcal{V}(\mathcal{M}(3))=\mathcal{V}$. \ In this case, there is a
column relation in $\mathcal{C}_{\mathcal{M}(4)}$ of the form 
\begin{equation*}
Y^{2}X^{2}=\alpha _{1}1+\alpha _{2}X+\alpha _{3}Y+\alpha _{4}X^{2}+\alpha
_{5}YX+\alpha _{6}Y^{2}+\alpha _{7}YX^{2}+\alpha _{8}Y^{2}X.
\end{equation*}%
To compute $\alpha _{1},\cdots ,\alpha _{8}$, let $\mathbf{v}$ denote the
compression of column $Y^{2}X^{2}$ in $\mathcal{M}(4)$ to rows indexed by
the basis $\mathcal{B}$, i.e., $\mathbf{v:}=(\beta _{22},\beta _{32},\beta
_{23},\beta _{42},\beta _{33},\beta _{24},\beta _{43},\beta _{34})^{T}$. \
Since $\mathcal{M}(4)$ is recursively generated, we have $X^{4}=YX$ and $%
YX^{3}=Y^{2}$ in $\mathcal{C}_{\mathcal{M}(4)}$, whence $\beta _{43}=\beta
_{14}$ and $\beta _{34}=\beta _{05}$. \ Thus $\mathbf{v}$ is expressed in
terms of the original data from $\beta ^{(6)}$. \ Let $J$ denote the
compression of $\mathcal{M}(3)$ to rows and columns indexed by elements of $%
\mathcal{B}$; then $J$ is invertible and $\alpha :=(\alpha _{1},\cdots
,\alpha _{8})$ is uniquely determined by 
\begin{equation}
\alpha ^{T}=J^{-1}\mathbf{v}.  \label{eq67new}
\end{equation}%
Now let 
\begin{equation*}
k(x,y):=y^{2}x^{2}-(\alpha _{1}+\alpha _{2}x+\alpha _{3}y+\alpha
_{4}x^{2}+\alpha _{5}yx+\alpha _{6}y^{2}+\alpha _{7}yx^{2}+\alpha _{8}y^{2}x.
\end{equation*}%
Since $k(X,Y)=0$ in $\mathcal{C}_{\mathcal{M}(4)}$, then $k|_{\mathcal{V}%
}\equiv 0$, so it follows from (\ref{eq62new}) and (\ref{eq62prime}) that $%
k\equiv h$, whence $\Lambda _{\beta }(k)=0$. \ Thus, if $k$ is computed as
above (using (\ref{eq67new})) and $\Lambda _{\beta }(k)\neq 0$, then $\beta $
has no representing measure.\newline
(ii) \ Let $k(x,y)$ be computed as above. \ Even without knowing the points
of $\mathcal{V}$ explicitly, if we know that $k|_{\mathcal{V}}\equiv 0$,
then from Lemma \ref{lem64}, Lemma \ref{lem31}, and (\ref{eq62prime}) it
follows that $k=h$, so $\beta $ has representing measure if and only if $%
\Lambda _{\beta }(k)=0$.\newline
(iii) \ Finally, we note that for the extremal problem for $\mathcal{M}(3)$
with $Y=X^{3}$, $\mathcal{M}(3)\geq 0$, $\mathcal{M}(2)>0$ and $r=v=8$, we
can always assume that $\mathcal{B}:=\{1,X,Y,X^{2},YX,Y^{2},$ $%
YX^{2},Y^{2}X\}$ is a basis for $\mathcal{C}_{\mathcal{M}(3)}$. \ Indeed,
suppose that a maximal linearly independent set of columns is $%
\{1,X,Y,X^{2},YX,Y^{2},YX^{2},Y^{3}\}$. \ Then there is a column relation of
the form $Y^{2}X=\alpha _{1}YX^{2}+\alpha _{2}Y^{3}+p(X,Y)\;(\deg \;p\leq 2)$%
. $\ $If $\alpha _{2}=0$, then (since $Y=X^{3})$, $\mathcal{V}(\mathcal{M}%
(3))$ is a subset of the zeros of $x^{7}=\alpha _{1}x^{5}+p(x,x^{3})$,
whence $v\leq 7$, a contradiction. \ Thus, $\alpha _{2}\neq 0$, and since $%
r=8$, it follows that\ $\mathcal{B}$ is a basis. \ A similar argument can be
used in the case when $\{1,X,Y,$ $X^{2},YX,Y^{2},$ $Y^{2}X,Y^{3}\}$ is a
basis. \ This completes the analysis of the extremal problem (\ref{eq61new}).
\end{remark}

\section{\label{sectexample}An Example with $r<v<+\infty $}

In this section we present an example in which we solve a truncated moment
problem with $r<v<+\infty $. Based on a number of examples and results in %
\cite{tcmp3}, \cite{tcmp7} and \cite{tcmp9}, we conjecture that in such
cases, if $\mathcal{M}(n)(\beta )$ has a representing measure, then a
minimal representing measure is $v$-atomic, and corresponds to a rank-$v$
positive extension $\mathcal{M}(n+k)$ (for some $k\leq v-r$), followed by a
flat extension $\mathcal{M}(n+k+1)$. \ In \cite{Fia4} we present an
algorithm for determining the existence of representing measures in a broad
class of truncated moment problems with $r<v<+\infty $; the following
example may be viewed as an instance of this algorithm, and also illustrates
Proposition \ref{cor36}.

\begin{example}
\label{ex71}Consider 
\begin{equation*}
{\mathcal{M}}(3)=\left( 
\begin{array}{cccccccccc}
1 & 0 & 0 & 1 & 2 & 5 & 0 & 0 & 0 & 0 \\ 
0 & 1 & 2 & 0 & 0 & 0 & 2 & 5 & 14 & 42 \\ 
0 & 2 & 5 & 0 & 0 & 0 & 5 & 14 & 42 & 132 \\ 
1 & 0 & 0 & 2 & 5 & 14 & 0 & 0 & 0 & 0 \\ 
2 & 0 & 0 & 5 & 14 & 42 & 0 & 0 & 0 & 0 \\ 
5 & 0 & 0 & 14 & 42 & 132 & 0 & 0 & 0 & 0 \\ 
0 & 2 & 5 & 0 & 0 & 0 & 5 & 14 & 42 & 132 \\ 
0 & 5 & 14 & 0 & 0 & 0 & 14 & 42 & 132 & 429 \\ 
0 & 14 & 42 & 0 & 0 & 0 & 42 & 132 & 429 & 2000 \\ 
0 & 42 & 132 & 0 & 0 & 0 & 132 & 429 & 2000 & 338881%
\end{array}%
\right) .
\end{equation*}%
We have ${\mathcal{M}}(3)\geq 0$, ${\mathcal{M}}(2)>0$ (positive and
invertible), and $r:=\operatorname{rank}\;{\mathcal{M}}(3)=8$, with column
relations 
\begin{equation}
Y=X^{3},  \label{eq61}
\end{equation}%
and 
\begin{equation}
Y^{3}=q(X,Y),  \label{eq62}
\end{equation}%
where $q(x,y):=-2285x+5720y-34441yx^{2}+578y^{2}x$. \ Let $%
r_{1}(x,y):=y-x^{3}$ and $r_{2}(x,y):=y^{3}+2285x-5720y+3441x^{2}y-578xy^{2}$%
. \ Then $\ker {\mathcal{M}}(3)=\langle \hat{r}_{1},~\hat{r}_{2}\rangle $
and ${\mathcal{V}}_{\beta }\equiv \{(x,y)\in \mathbb{R}%
^{2}:r_{1}(x,y)=r_{2}(x,y)=0\}$. \ A calculation shows that $v:=\operatorname{card%
}\;{\mathcal{V}}_{\beta }=9$. \ Now ${\mathcal{M}}(3)$ is positive,
recursively generated (trivially, because ${\mathcal{M}}(2)$ is invertible),
and $r<v$; further, Proposition \ref{cor36} implies that $\beta ^{(6)}$ is
consistent. \ We will show that the minimal representing measure for $\beta
^{(6)}$ is $v$-atomic (cf. Question \ref{quest1}).\newline
\indent%
If $\mu $ is a finitely atomic representing measure for $\beta $, then ${%
\mathcal{M}}(4)[\mu ]$ is recursively generated \cite{tcmp3}. \ Conversely,
any recursively generated extension ${\mathcal{M}}(4)$ of ${\mathcal{M}}(3)$
must satisfy 
\begin{equation}
X^{4}=YX  \label{eq63}
\end{equation}%
and 
\begin{equation}
YX^{3}=Y^{2}.  \label{eq64}
\end{equation}%
Further, since $Y^{3}=q(X,Y)$, in the column space of ${\mathcal{M}}(4)$ we
must have 
\begin{equation}
Y^{3}X=(xq)(X,Y)=q(X,Y)X  \label{eq65}
\end{equation}%
and 
\begin{equation}
Y^{4}=(yq)(X,Y)=Yq(X,Y).  \label{eq66}
\end{equation}%
\indent%
Using these column relations, we see that ${\mathcal{M}}(4)$ is completely
defined (i.e., all moments of degrees $7$ and $8$ are determined). \ On the
other hand, a calculation shows that in ${\mathcal{C}}_{{\mathcal{M}}(4)}$, $%
Y^{2}X^{2}$ is independent of ${\mathcal{B}}:=\{1,~X,~Y,~X^{2},~XY,~Y^{2},$ $%
X^{2}Y,~XY^{2}\}$. \ Thus, $\operatorname{rank}\;{\mathcal{M}}(4)=9$. \ Since a
flat extension of a positive, recursively generated moment matrix is
necessarily recursively generated \cite{tcmp3}, it follows that there is no
flat extension ${\mathcal{M}}(4)$ of ${\mathcal{M}}(3)$, and thus there is
no $8$-atomic representing measure for $\beta $.\newline
\indent%
Note that ${\mathcal{M}}(4)$ is extremal; indeed, the variety of ${\mathcal{M%
}}(4)$ consists of the common zeros of $r_{1},~r_{2},~xr_{1},$ $%
yr_{1},~xr_{2}$, and $yr_{2}$, and thus coincides with ${\mathcal{V}}_{\beta
}$ (which has $9$ points). \ Rather than using Theorem \ref{thm32}, we will
show that ${\mathcal{M}}(4)$ has a unique, $9$-atomic, representing measure
by a direct construction. \ Observe that relations (\ref{eq63})-(\ref{eq66}%
), together with recursiveness, completely determine any recursively
generated extension ${\mathcal{M}}(5)$ via the following relations: $%
X^{5}=YX^{2}$, $YX^{4}=Y^{2}X$, $Y^{2}X^{3}=Y^{3}$, $%
Y^{3}X^{2}=(x^{2}q)(X,Y) $, $Y^{4}X=(yxq)(X,Y)$, $Y^{5}=(y^{2}q)(X,Y)$. \ A
calculation of the degree-$5$ columns using these relations shows that these
columns do fit together to form a moment matrix, which is clearly a flat
(i.e., rank-preserving) extension of ${\mathcal{M}}(4)$. \ It thus follows
from \cite[Corollary 5.14]{tcmp1} that ${\mathcal{M}}(5)$ has a unique
representing measure $\mu $, which is $9$-atomic. \ Further, from the
recursive definition of ${\mathcal{M}}(5)$, it follows that $\mu $ is the
unique representing measure for $\beta $. \ \textbf{\qed}
\end{example}


\begin{thebibliography}{CuFi10}
\bibitem[AhKr]{AhKr} N.I. Ahiezer and M. Krein, \textit{Some Questions in
the Theory of Moments}, Transl. Math. Monographs, vol. 2, American Math.
Soc., Providence, 1962.

\bibitem[Akh]{Akh} N.I. Akhiezer, \textit{The Classical Moment Problem},
Hafner Publ. Co., New York, 1965.

\bibitem[Atk]{Atk} K. Atkinson, \textit{Introduction to Numerical Analysis},
Wiley and Sons, 2nd. Ed. \ 1989.

\bibitem[BaTe]{BaTe} C. Bayer and J. Teichmann, The proof of Tchakaloff's
Theorem, preprint 2005.

\bibitem[CLO]{CLO'S} D. Cox, J. Little and D. O'Shea, \textit{Ideals,
Varieties and Algorithms: An Introduction to Computational Algebraic
Geometry and Commutative Algebra}, Second Edition, Springer-Verlag, New
York, 1992.

\bibitem[CuFi1]{RGWSI} R. Curto and L. Fialkow, Recursively generated
weighted shifts and the subnormal completion problem, Integral Equations
Operator Theory 17(1993), 202-246.

\bibitem[CuFi2]{tcmp1} R. Curto and L. Fialkow, Solution of the truncated
complex moment problem with flat data, \textit{Memoirs Amer. Math. Soc}. no.
568, Amer. Math. Soc., Providence, 1996.

\bibitem[CuFi3]{tcmp2} R. Curto and L. Fialkow, Flat extensions of positive
moment matrices: Relations in analytic or conjugate terms, \textit{Operator
Th.: Adv. Appl}. 104(1998), 59-82.

\bibitem[CuFi4]{tcmp3} R. Curto and L. Fialkow, Flat extensions of positive
moment matrices: Recursively generated relations, \textit{Memoirs Amer.
Math. Soc}. no. 648, Amer. Math. Soc., Providence, 1998.

\bibitem[CuFi5]{tcmp4} R. Curto and L. Fialkow, The truncated complex $K$%
-moment problem, \textit{Trans. Amer. Math. Soc}. 352(2000), 2825-2855.

\bibitem[CuFi6]{tcmp5} R. Curto ad L. Fialkow, The quadratic moment problem
for the unit disk and unit circle, \textit{Integral Equations Operator Theory%
} 38(2000), 377-409.

\bibitem[CuFi7]{tcmp6} R. Curto and L. Fialkow, Solution of the singular
quartic moment problem, \textit{J. Operator Theory }48(2002), 315-354.

\bibitem[CuFi8]{tcmp7} R. Curto and L. Fialkow, Solution of the truncated
parabolic moment problem, \textit{Integral Equations Operator Theory }%
50(2004), 169-196.

\bibitem[CuFi9]{tcmp8} R. Curto and L. Fialkow, A duality proof of
Tchakaloff's theorem, \textit{J. Math. Anal. Appl}. 269(2002), 519-532.

\bibitem[CuFi10]{tcmp9} R. Curto and L. Fialkow, Solution of the truncated
hyperbolic moment problem, \textit{Integral Equations Operator Theory }%
52(2005), 181-218.

\bibitem[CuFi11]{tcmp10} R. Curto and L. Fialkow, Truncated $K$-moment
problems in several variables, \textit{J. Operator Theory 54(2005), 189-226}.

\bibitem[Dou]{Dou} R.G. Douglas, On majorization, factorization, and range
inclusion of operators on Hilbert spaces, \textit{Proc. Amer. Math. Soc}.
17(1966), 413-415.

\bibitem[Fia1]{FiCM} L. Fialkow, Positivity, extensions and the truncated
complex moment problem, \textit{Contemporary Math}. 185(1995), 133-150.

\bibitem[Fia2]{FiaOT} L. Fialkow, Minimal representing measures arising from
rank-increasing moment matrix extensions, \textit{J. Operator Theory}
42(1999), 425-436.

\bibitem[Fia3]{FiaIEOT} L. Fialkow, Truncated complex moment problems with a 
$Z\bar{Z}$ relation, \textit{Integral Equations Operator Theory} 45(2003),
405-435.

\bibitem[Fia4]{Fia4} L. Fialkow, Extensions and varieties of moment
matrices, in preparation.

\bibitem[FiPe]{FP} L. Fialkow and S. Petrovic, A moment matrix approach to
multivariable cubature, \textit{Integral Equations Operator Theory }%
52(2005), 85-124.

\bibitem[Ful]{Ful} W. Fulton, \textit{Algebraic curves. An introduction to
algebraic geometry, }Mathematics Lecture Notes Series. W. A. Benjamin, Inc.,
New York-Amsterdam, 1969.

\bibitem[KrNu]{KrNu} M.G. Krein and A.A. Nudel'man, \textit{The Markov
Moment Problem and Extremal Problems}, \textit{Transl. Math. Monographs, }%
vol. 50, American Mathematical Society, Providence, R.I., 1977.

\bibitem[KuMa]{KuM} S. Kuhlmann and M. Marshall, Positivity, sums of squares
and the multidimensional moment problem, \textit{Trans. Amer. Math. Soc}.
354(2002), 4285-4301.

\bibitem[Lau1]{Lau1} M. Laurent, Semidefinite representations for finite
varieties, preprint 2002.

\bibitem[Lau2]{Lau2} M. Laurent, Revisiting two theorems of Curto and
Fialkow on moment matrices, \textit{Proc. Amer. Math. Soc}. 133 (2005),
2965-2976.

\bibitem[MMM]{MMM} M.G. Marinari, H.M. M\"{o}ller, and T. Mora, On
multiplicities in polynomial system solving, \textit{Trans. Amer. Math. Soc.}
348(1996), 3283 -- 3321.

\bibitem[Moe1]{Moe} H.M. M\"{o}ller, An inverse problem for cubature
formulae, \textit{Computat. Technol.} 9 (2004), 13 -- 20.

\bibitem[Moe2]{Moe2} H.M. M\"{o}ller, On square positive extensions and
cubature formulas, to appear in \textit{J. Comput. Applied Math}., Special
issue ed. by W. zu Castell and F. Filbir.

\bibitem[MoSa]{MoSa} H.M. M\"{o}ller and T. Sauer, $H$-bases for polynomial
interpolation and system solving. Multivariate polynomial interpolation. 
\textit{Adv. Comput. Math}. 12 (2000), 335--362.

\bibitem[Mys]{Mys} I.P. Mysovskikh, On Chakalov's theorem, \textit{U.S.S.R.
Comput.\ Math.\ and Math.\ Phys.}\ 15 (1975), 221--227; (translation in 
\textit{Zh.\ Vychisl.\ Mat.\ i Mat.\ Fiz}.\ 15 (1975), 1589--1593).

\bibitem[PoSc]{PoS} V. Powers and C. Scheiderer, The moment problem for
non-compact semialgebraic sets, \textit{Advances in Geometry} 1(2001), 71-88.

\bibitem[Put]{tchaka} M. Putinar, On Tchakaloff's Theorem, \textit{Proc.
Amer. Math. Soc}. 125(1997), 2409-2414.

\bibitem[PuVa]{PuVa} M. Putinar and F.-H. Vasilescu, Solving moment problems
by dimensional extension, \textit{Ann. of Math}. (2) 149(1999), no. 3,
1087-1107.

\bibitem[Sche]{Sche} C. Scheiderer, Sums of squares of regular functions on
real algebraic varieties, \textit{Trans. Amer. Math. Soc}. 352(2000),
1039-1069.

\bibitem[Schm1]{Sch} K. Schm\"{u}dgen, The $K$-moment problem for
semi-algebraic sets, \textit{Math. Ann}. 289(1991), 203-206.

\bibitem[Schm2]{SchNew} K. Schm\"{u}dgen, On the moment problem of closed
semi-algebraic sets, \textit{J. Reine Angew. Math}. 558 (2003), 225--234.

\bibitem[ShTa]{ShTa} J.A. Shohat and J.D. Tamarkin, \textit{The Problem of
Moments}, Math. Surveys I, American Math. Soc., Providence, 1943.

\bibitem[Smu]{Smu} J.L. Smul'jan, An operator Hellinger integral (Russian), 
\textit{Mat. Sb}. 91(1959), 381-430.

\bibitem[Sto1]{Sto1} J. Stochel, Moment functions on real algebraic sets, 
\textit{Ark. Mat}. 30(1992), 133-148.

\bibitem[Sto2]{Sto2} J. Stochel, Solving the truncated moment problem solves
the moment problem, \textit{Glasgow J. Math}. 43(2001), 335-341.

\bibitem[StSz1]{StSz1} J. Stochel and F.H. Szafraniec, Algebraic operators
and moments on algebraic sets, \textit{Portugal. Math.} 51(1994), 25-45.

\bibitem[StSz2]{StSz2} J. Stochel and F.H. Szafraniec, The complex moment
problem and subnormality: A polar decomposition approach, \textit{J. Funct.
Anal}. 159(1998), 432-491.

\bibitem[Tch]{Tch} V. Tchakaloff, Formules de cubatures m\'{e}caniques \`{a}
coefficients non n\'{e}gatifs, \textit{Bull. Sc. Math}. 81(1957).

\bibitem[Wol]{Wol} Wolfram Research, Inc., \textit{Mathematica}, Version
4.2, Wolfram Research, Inc., Champaign, IL, 2002.
\end{thebibliography}
\end{document}